\documentclass[12pt,reqno]{amsart}
\usepackage{mathrsfs}
\usepackage{graphicx,epsfig,amsmath,amsfonts,dsfont}
\usepackage{color}
\usepackage{latexsym,amssymb}
\usepackage{lineno,version}
\usepackage{multirow}
\usepackage{bm}
\usepackage{enumitem}
\usepackage{float}
\usepackage{booktabs}
\usepackage{hyperref} 
\usepackage{geometry}
\allowdisplaybreaks

\catcode`\@=11
\theoremstyle{plain}
\@addtoreset{equation}{section}   

\@addtoreset{figure}{section}
\renewcommand\thefigure{\thesection.\@arabic\c@figure}

\newtheorem{theorem}{\noindent Theorem}[section]
\newtheorem{remark}{\noindent Remark}[section]
\newtheorem{algorithm}{\noindent Algorithm}[section]
\newtheorem{proposition}{\noindent Proposition}[section]
\newtheorem{lemma}{\noindent Lemma}[section]
\newtheorem{example}{\noindent Example}[section]
\newtheorem{definition}{\noindent Definition}[section]

\newtheorem*{prob}{\noindent IP}

\newcommand{\ba}{\begin{array}}\newcommand{\ea}{\end{array}}
\newcommand{\be}{\begin{eqnarray}}\newcommand{\ee}{\end{eqnarray}}
\newcommand{\beq}{\begin{equation}}\newcommand{\eeq}{\end{equation}}
\newcommand{\bex}{\begin{eqnarray*}}
\newcommand{\eex}{\end{eqnarray*}}

\font\tenbi=cmmib10   at 11 pt
\font\sevenbi=cmmib10 at 9pt
\font\fivebi=cmmib7 at 6pt
\newfam\bifam
\textfont\bifam=\tenbi \scriptfont\bifam=\sevenbi
\scriptscriptfont\bifam=\fivebi

\font\tendb=msbm10 at 12 pt
\font\sevendb=msbm7
\newfam\dbfam
\textfont\dbfam=\tendb \scriptfont\dbfam=\sevendb

\def\div{\mbox{div}}

\begin{document}

\title[]
{Total variation regularization for recovering the spatial source term in a time-fractional diffusion equation$^*$}
\author[]
{Bin Fan$^{1,\dag}$}
\thanks{\hskip -12pt
${}^*$This work was supported by the Natural Science Foundation of Fujian Province (Grant No. 2022J05198), the Foundation of Fujian Provincial Department of Education (Grant No. JAT210293), and the Natural Science Foundation of Fujian University of Technology (Grant No. GY-Z21069).\\
${}^{1}$School of Computer Science and Mathematics, Fujian University of Technology, 350118 Fuzhou, China.\\
Email: bfan@fjut.edu.cn}



\maketitle

\begin{abstract}
In this paper, we consider an inverse space-dependent source problem for a time-fractional diffusion equation.
To deal with the ill-posedness of the problem, we transform the problem into  an optimal control problem with total variational (TV) regularization.
In contrast to the classical Tikhonov model incorporating $L^2$ penalty terms, the inclusion of a TV term proves advantageous in reconstructing solutions that exhibit discontinuities or piecewise constancy.
The control problem is approximated by a fully discrete scheme, and convergence results are provided within this framework.
Furthermore, a linearized primal-dual iterative algorithm is proposed to solve the discrete control model based on an equivalent saddle-point reformulation, and several numerical experiments are presented to demonstrate the efficiency of the algorithm.

\end{abstract}

\section{Introduction}\label{sec:1}

Let $T>0$ and $\Omega\subset\mathbb{R}^d$ ($d=1,2,3$) be a bounded convex polygonal domain. In this paper, we are interested in the following initial-boundary value problem for a time-fractional diffusion equation (TFDE) with the homogeneous Neumann boundary condition
\begin{eqnarray}
   && \left\{\begin{array}{ll}
               \partial_t^\alpha u(x,t)+L u(x,t)=\mu(t)f(x), & (x,t)\in D_T:=\Omega\times (0,T), \\
               u(x,0)=0,& x\in\Omega, \\
               \partial_\nu u(x,t)=0, &   (x,t)\in\partial\Omega\times(0,T),
             \end{array}\right. \label{general-dir-prob}
\end{eqnarray}
where
\begin{itemize}[]
\item[$\circ$]$\partial_t^\alpha$ denotes the Caputo fractional derivative of order $\alpha$ defined by (see, e.g., \cite{podlubny1998fractional,diethelm2002analysis})
\begin{eqnarray}
 \partial_t^\alpha u(x,t):=\frac{1}{\Gamma(1-\alpha)}\int_0^t \frac{\partial u(x,s)}{\partial s} \frac{\mathrm{d}s}{(t-s)^{\alpha}},\quad 0<\alpha<1,\quad 0<t\leq T
\label{det-caputo-1}
\end{eqnarray}
with $\Gamma(\cdot)$ being the Gamma function.
\item[$\circ$]$L$ is a symmetric uniformly elliptic operator defined on $H^2(\Omega)$ given by\begin{eqnarray*}
 Lv(x,t):= - \sum\limits_{i=1}^d\frac{\partial}{\partial x_i}\left(\sum\limits_{i=1}^da_{ij}(x)\frac{\partial}{\partial x_j}v(x,t) \right)+b(x)u(x,t),
\end{eqnarray*}
in which the coefficients satisfy
\begin{eqnarray*}
   &&  a_{ij}=a_{ji},\quad a_{ij}\in C^1(\bar{\Omega}),\quad 1\leq i,j\leq d,\\
   && a_0\sum\limits_{i=1}^d\xi_i^2\leq \sum\limits_{i,j=1}^da_{ij}(x)\xi_i\xi_j,\quad \forall x\in\bar{\Omega},
   \quad \forall \xi\in\mathds{R}^d,\quad a_0>0,\\
   && b\in C(\bar{\Omega}),\quad b(x)\ge 0,\quad \forall x\in\bar{\Omega}.
\end{eqnarray*}
\item[$\circ$]$\partial_\nu$ is the unit outward normal derivative on the boundary $\partial\Omega$.
\end{itemize}
The source term in equation \eqref{general-dir-prob} is expressed as separated variables, where the spatial component $f(x)$ represents the spatial distribution, such as that of the contaminant source, and the temporal component $\mu(t)$ characterizes the temporal evolution pattern.
The governing equation presented in  \eqref{general-dir-prob} is a prototypical example of a broad class of TFDEs, which have been proposed as a potent tool for describing anomalous diffusion phenomena in heterogeneous media.
In recent decades, TFDEs have been extensively investigated from both theoretical and numerical perspectives; see, e.g., recent papers \cite{jin2016two,jin2019numerical,zhu2019fast} and the references therein.

However, for some practical problems, the initial data, diffusion coefficients, or source
term may be unknown and we want to determine them by additional measurement data which will yield to
some fractional diffusion inverse problems.
In \cite{sakamoto2011initial}, Sakamoto and  Yamamoto established several uniqueness results for various inverse problems of TFDE.
As we know, these problems are generally ill-posed, so the so-called regularization technique is necessary to solve the inverse problem.
By utilizing the spectral information of the operator $-L$, Wei et al. \cite{wei2016inverse} established both uniqueness and a stability estimate for solving the inverse time-dependent source problem from boundary Cauchy data. They employed the Tikhonov regularization method to tackle the inverse source problem and proposed a conjugate gradient algorithm to obtain an excellent approximation to the minimizer of the Tikhonov regularization functional.
Jiang et al. \cite{jiang2017weak} established the uniqueness of the inverse space-dependent source problem from partial interior observations and proposed an iterative thresholding algorithm based on Tikhonov regularization. Subsequently, a detailed  analysis of the numerical restoration algorithm was presented in \cite{jiang2020numerical}.
Based on the variational methods, several optimal control models with Tikhonov regularization have been proposed and analyzed in recent years for solving TFDEs-based inverse problems, see e.g., \cite{ye2013spectral,zhou2016finite,jin2020pointwise,hrizi2023reconstruction}.
For a more comprehensive bibliography on inverse problems for TFDEs, we refer e.g., to \cite{jin2015tutorial,liu2019inverse}.

The classical Tikhonov methods with standard penalty terms, such as $\|\cdot\|_{L^2(\Omega)}^2$, are known to suffer from the over-smoothing solutions. To prevent over-smoothing of the solution by the Tikhonov method, there are two popular approaches: one is the least square method with total variation regularization terms \cite{chavent1997regularization,wang2013total,Vasin2016Regularization}; the other is the novel fractional regularization methods \cite{2008Regularization,2011Fractional}.
The latter is  typically predicated on spectral decomposition of solution operators, otherwise it is difficult to have a clear explanation of fractional powers.
TV regularization is highly advantageous for recovering non-smooth or discontinuous solutions and has wide applications in the field of image processing.
However, its application in the inverse problem of differential equations remains limited, see, e.g., \cite{wang2013total,hinze2018identifying,hinze2019finite,fan2021identifying}.
One of the reasons for this limitation is attributed to the non-differentiability of the TV-term, which poses challenges for both theoretical analysis and numerical computation.

In this paper, our primary focus is on recovering a potentially non-smooth source function $f(x)$ from partial observations. While numerous studies have addressed source identification using observations spanning the entire domain or at the final moment, there are relatively few studies on partial observations or the recovery of non-smooth solutions.
To achieve this, we will make attempt to construct an optimal control problem with TV regularization for solving the inverse source problem.
Then the control problem is discretized by using the standard piecewise linear finite element in space and a finite difference scheme in time.
Unlike many of the previous literatures, we establish some new convergence results in the finite dimensional space, see, Theorems \ref{thm-convfm-ipthm} and \ref{thm-discon}.
Furthermore, we present an effective restoration algorithm and report the corresponding numerical results.

The rest of this paper is organized as follows.
In Section \ref{sec:2}, we formally present the inverse source problem of TFDE, provide some notations for this paper, and list some useful auxiliary results. In Section \ref{sec:3}, the inverse source problem is first reformulated as an optimal control problem with TV regularization, and the existence and stability results are provided.
In Section \ref{sec:4}, the discrete form of the optimal control problem is obtained by employing the standard Galerkin method with finite element in space and the finite difference scheme in time. Then we establish some convergence results for the proposed model under appropriate assumptions.
In Section \ref{sec:5}, a linearized primal-dual algorithm is proposed to effectively solve the discrete constrained minimization problem.
Numerical implementation and some results are reported in Section \ref{sec:6}. Finally, we give some concluding remarks in the last section.

\section{Preliminaries}\label{sec:2}
Without loss of generality, we exclusively consider the following equation due to its simplicity in numerical simulation and our conviction that the underlying ill-posedness is fundamentally identical.
\begin{equation}\label{prob-forward}
  \left\{\begin{array}{ll}
           (\partial_t^{\alpha}-\Delta+1)u(x,t)=\mu(t)f(x), & (x,t)\in \Omega\times (0,T), \\
           u(x,0)=0, & x\in\Omega,  \\
           \partial_\nu u(x,t)=0, & (x,t)\in\partial\Omega\times(0,T).
         \end{array}\right.
\end{equation}
Throughout the paper, we make the assumption that the function $\mu$ is a known function on $[0,T]$. Additionally, we denote $u(f)$ as the solution of equation \eqref{prob-forward} with respect to space source term $f$.
The primary focus of this paper lies in the numerical aspect of the following inverse source problem.
\begin{prob}
Let $\omega\subset\Omega$ be a nonempty open subdomain, determine the spatial component of source term $f$ in \eqref{prob-forward} from a partial interior measurement $u^\delta\in L^2(\omega\times(0,T))$ of the exact data $u^\dag:=u(f^*)|_{\omega}$ satisfying
\begin{equation}\label{noisy-data}
  \|u^\delta-u^\dag\|_{ L^2(\omega\times(0,T))}\leq\delta,
\end{equation}
where $f^*\in L^2(\Omega)$ and $\delta>0$ stand for the true source term and the noise level, respectively.
\end{prob}

The inverse problem we are considering is ill-posed. This aspect
can be well explained by spectral decomposition of the solution operator $u(f)$, for example, the spectral decomposition of solution operator typically
encompasses Mittag-Leffler function terms that demonstrate significant decay, thereby leading
to potential instability in direct solutions. For more detail, we highly recommend consulting
references [1], [2], among others.

For the solution regularity of the direct problem \eqref{prob-forward}, we define the fractional Sobolev spaces $_0\!H^\alpha (0,T)$ as (see, e.g., \cite{kubica2020time})
\begin{equation*}
  _0\!H^\alpha (0,T):=\left\{\begin{array}{ll}
                               H^\alpha(0,T), & 0<\alpha<1/2, \\[2pt]
                               \left\{g\in H^{1/2}(0,T)~|~\int_0^T\frac{|g(t)|^2}{t}\mathrm{d}t<\infty\right\}, & \alpha=1/2, \\[2pt]
                               \left\{g\in H^\alpha(0,T)~|~g(0)=0\right\}, & 1/2<\alpha<1.
                             \end{array}\right.
\end{equation*}
Then we have the following regularity results for the direct problem.
\begin{lemma}\label{lem-dirpb}$($\cite[Lemma 2.4]{jiang2017weak} and \cite[Lemma 2.2]{jiang2020numerical}$)$
Let $f\in L^2(\Omega)$ and $\mu\in L^\infty(0,T)$. Then the initial-boundary value problem  \eqref{prob-forward} admits a unique solution $u(f)\in _0\!H^\alpha (0,T;L^2(\Omega))\cap C([0,T];L^2(\Omega))\cap L^2(0,T;H^2(\Omega))$ such that
$$\|u(f)\|_{_0\!H^\alpha (0,T;L^2(\Omega))} + \|u(f)\|_{C([0,T];L^2(\Omega))} + \|u(f)\|_{L^2(0,T;H^2(\Omega))}\leq c\|f\|_{L^2(\Omega)},$$
where constant $c>0$ depending on $\Omega$, $T$, $\alpha$ and $\mu$.
\end{lemma}

\begin{definition}\label{def-frac-ID}
Let $\alpha\in(0,1)$ and $t\in(0,T]$. $ D_{T^-}^\alpha$ and $ I_{T^-}^\alpha$ are the right-sided Riemann-Liouville fractional derivative and integral, respectively, defined by
\begin{equation*}
   D_{T^-}^\alpha u(x,t):=-\frac{1}{\Gamma(1-\alpha)} \frac{\mathrm{d}}{\mathrm{d}t} \int_t^T \frac{u(x,s)}{(s-t)^\alpha}\mathrm{d}s, \qquad  I_{T^-}^\alpha u(x,t):=\frac{1}{\Gamma(\alpha)}  \int_t^T \frac{u(x,s)}{(s-t)^{1-\alpha}}\mathrm{d}s.
\end{equation*}
 $ D_{0^+}^\alpha$ and $ I_{0^+}^\alpha$ are the left-sided Riemann-Liouville fractional derivative and integral, respectively, defined by
\begin{equation*}
  D_{0^+}^\alpha u(x,t):=\frac{1}{\Gamma(1-\alpha)} \frac{\mathrm{d}}{\mathrm{d}t} \int_0^t \frac{u(x,s)}{(t-s)^\alpha}\mathrm{d}s,\qquad I_{0^+}^\alpha u(x,t):=\frac{1}{\Gamma(\alpha)}  \int_0^t \frac{u(x,s)}{(t-s)^{1-\alpha}}\mathrm{d}s.
\end{equation*}
\end{definition}
It is well known that the Riemann-Liouville and Caputo fractional derivatives have the following relationship
\begin{equation}\label{rel-RandC}
 D_{0^+}^\alpha u(x,t)=\partial_t^\alpha u(x,t) + \frac{u(x,0)t^{-\alpha}}{\Gamma(1-\alpha)}.
\end{equation}

Throughout the paper, we employ the symbol $c$, with or without subscripts or bars, to denote generic positive constants that may vary across different instances.
The symbols $\to$ and $\rightharpoonup$ are employed to denote the strong and weak convergence of sequences of functions in the respective space.
Let $L^p(\Omega)$ ($1\leq p<\infty$) be the classical Lebesgue spaces, and specifically, the commonly used $L^2$-space is equipped with the inner product $(\cdot,\cdot)$. We also use the standard notations of Sobolev spaces $H^2(\Omega)$, $H_0^1(\Omega)$ and $H^\alpha$ ($\alpha\in(0,1)$), with the corresponding norms (see, e.g., \cite{1975Sobolev}).

\section{Optimal control problem with TV regularization}\label{sec:3}
We begin by briefly revisiting the space of functions having bounded total variation.
The definition of total variation can be formulated as
\begin{eqnarray*}
TV(f):=\sup\limits_{p\in\mathcal{B}}\int_\Omega f\div p \mathrm{d}x<\infty \quad\text{for}\quad \mathcal{B}:= \left\{p~|~p\in C_c^1(\Omega)^d,~|p|_\infty\leq 1 \right\}
\end{eqnarray*}
where $|p|_\infty:=\sup\limits_{x\in \Omega}(\sum\limits_{i=1}^d|p_i(x)|^2)^{1/2}$,
$\div$ denotes the divergence operator, and $C_c^1(\Omega)$ is the space of
continuously differentiable $\mathds{R}^d$-valued functions with compact support in $\Omega$.
The space of all functions in $L^1(\Omega)$ with bounded variation is denoted by
\begin{eqnarray*}
BV(\Omega):=\{f\in L^1(\Omega)~|~TV(f)<\infty\}.
\end{eqnarray*}
It can be proved that $BV(\Omega)$ is a Banach space equipped with the norm $\|f\|_{BV(\Omega)}:=\|f\|_{L^1(\Omega)} + TV(f)$.

\begin{proposition}$($\cite[Proposition 1.1]{chavent1997regularization}$)$\label{prop-BVfun}

$(i)$~ $($Lower semicontinuity$)$ If $\{f_m\}_{m=1}^\infty\subset BV(\Omega)$ and $f_m\to f$ in $L^1(\Omega)$ as $m\to \infty$, then $f\in BV(\Omega)$ and
\begin{eqnarray*}
   &&  TV(f)\leq\liminf\limits_{m\to\infty}TV(f_m).
\end{eqnarray*}

$(ii)$~For every bounded sequence $\{f_m\}_{m=1}^\infty\subset BV(\Omega)$ there exists a subsequence $\{f_{m_k}\}_{k=1}^\infty$ and $f\in BV(\Omega)$ such that $f_{m_k}\to f$ in $L^p(\Omega)$ as $k\to\infty$, $p\in[1,\frac{d}{d-1})$, if $d\geq 2;$ and $f_{m_k}\to f$ in $L^p(\Omega)$ as $k\to\infty$, $p\in[1,\infty)$, if $d= 1$.
\end{proposition}

(\textbf{IP}) can be reformulated as the integral equation
\begin{equation}
  \int_0^T\int_{\omega}(u(f)-u^\delta)^2\mathrm{d}x\mathrm{d}t=\|u(f)-u^\delta\|_{L^2(\omega\times(0,T))}^2=0,\label{prob-unreint}
\end{equation}
where $u(f)$ satisfies the initial condition $u(f)(\cdot,0)=0$ and the weak variational system of \eqref{prob-forward}
\begin{equation*}
  \int_0^T\int_{\Omega}(\partial_t^\alpha u(f)\varphi+\nabla u(f)\cdot\nabla\varphi+u(f)\varphi)\mathrm{d}x\mathrm{d}t=\int_0^T\int_{\Omega}\mu f\varphi\mathrm{d}x\mathrm{d}t
\end{equation*}
for any $\varphi\in H^\alpha(0,T;L^2(\Omega))\cap L^2(0,T;H^1(\Omega))$. To deal with the ill-posedness of problem \eqref{prob-unreint}, we consider the following optimal control problem with TV regularization
\begin{equation*}
 \min\limits_{f\in F_{ad}} \mathcal{J}(f),\quad \mathcal{J}(f):=\frac{1}{2}\|u(f)-u^\delta\|_{ L^2(\omega\times(0,T))}^2  + \gamma TV(f),\tag{$\mathcal{P}$}
\end{equation*}
where $\gamma>0$ is the regularization parameter, and admissible domain
\begin{equation*}
  F_{ad}:=\{f\in BV(\Omega)~|~-\infty<\underline{f}\leq f(x)\leq \overline{f}<\infty\quad\text{for a.e. in~~}\Omega\},
\end{equation*}
with $\underline{f}$, $\overline{f}$ being given constants. Notice here that the  uniform  boundedness of $f$ guarantees that $f\in L^s(\Omega)$ for all $1\leq s<\infty$.

\begin{theorem}\label{thm-exist-cp}
The optimization problem $(\mathcal{P})$ admits a minimizer for any $u^\delta\in L^2(\omega\times(0,T))$.
\end{theorem}
\noindent \emph{Proof }
By virtue of the nonnegativity of $\mathcal{J}(f)$, it is evident that $\inf_{f\in F_{ad}}\mathcal{J}(f)$ is finite. Therefore, there exists a minimizing sequence $\{f_m\}_{m=1}^\infty\subset F_{ad}$ such that
\begin{equation*}
  \lim\limits_{m\to\infty}\mathcal{J}(f_m)=\inf\limits_{f\in F_{ad}}\mathcal{J}(f).
\end{equation*}
Based on the definition of the set of admissible solutions, we know that $\{f_m\}_{m=1}^\infty$ is uniformly bounded in $BV(\Omega)$.
By Proposition \ref{prop-BVfun}, there exists a (not relabeled) subsequence of $\{f_m\}_{m=1}^\infty$ and an element $f^*\in BV(\Omega)$ such that
\begin{equation}\label{thm-exist-cp-eq1}
  f_m\to f^*\quad \text{in}\ L^1(\Omega)\ \text{as }m\to\infty,
\end{equation}
and
\begin{equation}\label{thm-exist-cp-eq2}
   TV(f^*)\leq\liminf\limits_{m\to\infty}TV(f_m).
\end{equation}
Since $\underline{f}\leq f_m(x)\leq \overline{f}$ for all $m\in\mathbb{N}$ and a.e. in $\Omega$, thus we obtain $\underline{f}\leq f^*(x)\leq \overline{f}$ a.e. in $\Omega$ by sending $m$ to $\infty$, which implies that $f^*\in F_{ad}$.
Furthermore, we observe that \eqref{thm-exist-cp-eq1} yields the convergence in the $L^s(\Omega)$-norm for all $1\leq s<\infty$. In fact, we have the following estimate
\begin{equation*}
  \|f_m-f^*\|_{L^s(\Omega}^s=\int_{\Omega}|f_m-f^*|\cdot|f_m-f^*|^{s-1}\mathrm{d}x\leq 2\overline{f}^{s-1}\|f_m-f^*\|_{L^1(\Omega)}.
\end{equation*}
Next we want to show that $f^*$ is a minimizer of optimization problem $(\mathcal{P})$.
It follows from Lemma \ref{lem-dirpb} that  the sequence $\{u(f_m)\}_{m=1}^\infty$ is uniformly bounded in $_0\!H^\alpha (0,T;L^2(\Omega))\cap L^2(0,T;H^2(\Omega))$ with respect to $m$. This indicates the existence of some $u^*\in _0\!H^\alpha (0,T;L^2(\Omega))\cap L^2(0,T;H^2(\Omega))$ and a (not relabeled) subsequence of $\{u(f_m)\}_{m=1}^\infty$, such that
\begin{equation*}
  u(f_m)\rightharpoonup u^*\quad \text{in}\  _0\!H^\alpha (0,T;L^2(\Omega))\cap L^2(0,T;H^2(\Omega))\ \text{as }m\to\infty.
\end{equation*}
Then we immediately get $u^*=u(f^*)$ according to the proof of \cite[Theorem 2.1]{jiang2020numerical}.
Finally, we use  \eqref{thm-exist-cp-eq2} and the lower semi-continuity of the $L^2$-norm  to conclude
\begin{align*}
\mathcal{J}(f^*) &= \frac{1}{2}\|u(f^*)-u^\delta\|_{ L^2(\omega\times(0,T))}^2 +  \gamma TV(f^*) \\
                &\leq\frac{1}{2}\liminf\limits_{m\to\infty} \|u(f_m)-u^\delta\|_{ L^2(\omega\times(0,T))}^2  + \gamma \liminf\limits_{m\to\infty}TV(f_m) \\
                & =\liminf\limits_{m\to\infty} \mathcal{J}(f_m)=\inf\limits_{f\in F_{ad}}\mathcal{J}(f),
\end{align*}
thus $f^*$ is indeed a minimizer to the optimization problem $(\mathcal{P})$.
\hfill$\Box$

Next, we establish the stability of problem $(\mathcal{P})$, demonstrating that the minimization problem $(\mathcal{P})$ effectively stabilizes the considered $(\mathcal{P})$ against perturbations in observation data within $\omega\times(0,T)$.

\begin{theorem}\label{thm-stab-cp}
Let $\{u_l^\delta\}_{l=1}^\infty\subset L^2(\omega\times(0,T))$ be a sequence satisfies
\begin{equation}\label{thm-stab-cp-eq1}
  u_l^\delta\to u^\delta\quad \text{in}\ L^2(0,T;L^2(\Omega))\ \text{as }l\to\infty,
\end{equation}
and $f_l$ is an arbitrary minimizer of the problem
\begin{equation*}
  \min\limits_{f\in F_{ad}} \mathcal{J}_l(f),\quad \mathcal{J}_l(f):=\frac{1}{2}\|u(f)-u_l^\delta\|_{ L^2(\omega\times(0,T))}^2  + \gamma TV(f)
\end{equation*}
for each $l\in\mathbb{N}$.
Then an unrelabeled subsequence of $\{f_l\}_{l=1}^\infty$ exists such that
\begin{equation}\label{thm-stab-cp-eq2}
  \lim\limits_{l\to\infty}\|f_l-f^*\|_{L^1(\Omega)}=0,\quad\text{and}\quad \lim\limits_{l\to\infty}TV(f_l)=TV(f^*),
\end{equation}
where $f^*\in F_{ad}$ is a minimizer of the optimization problem $(\mathcal{P})$.
\end{theorem}
\noindent \emph{Proof }
The existence of sequence $\{f_l\}_{l=1}^\infty\subset F_{ad}$ is guaranteed by Theorem \ref{thm-exist-cp}.
Then the sequence $\{f_l\}_{l=1}^\infty$ is uniformly bounded in $BV(\Omega)$ due to the optimality of $f_l$.
Reiterating the same argument as that presented in the proof of Theorem \ref{thm-exist-cp}, there exists a (not relabeled) subsequence of $\{f_l\}_{l=1}^\infty$ an element $f^*\in F_{ad}$ such that $f_l\to f^*$ in $L^1(\Omega)$ as $l\to\infty$ and  $TV(f^*)\leq\liminf_{l\to\infty}TV(f_l)$. Moreover, it can be inferred that
\begin{equation*}
  u(f_l)\rightharpoonup u(f^*)\quad \text{in}\  _0\!H^\alpha (0,T;L^2(\Omega))\cap L^2(0,T;H^2(\Omega))\ \text{as }l\to\infty.
\end{equation*}
up to taking a further (not relabeled) subsequence of $\{u(f_l)\}_{l=1}^\infty$.
In conjunction with \eqref{thm-stab-cp-eq1}, we can conclude that
\begin{equation*}
  u(f_l)-u_l^\delta\rightharpoonup u(f^*)-u^\delta\quad \text{in}\   L^2(0,T;L^2(\omega))\ \text{as }l\to\infty.
\end{equation*}
Consequently, for any $f\in F_{ad}$, again by the lower semi-continuity of the $L^2$-norm and TV term to deduce
\begin{align*}
\mathcal{J}(f^*) &= \frac{1}{2}\|u(f^*)-u^\delta\|_{ L^2(\omega\times(0,T))}^2 +  \gamma TV(f^*) \\
                &\leq\frac{1}{2}\liminf\limits_{l\to\infty} \|u(f_l)-u_l^\delta\|_{ L^2(\omega\times(0,T))}^2  + \gamma \liminf\limits_{l\to\infty}TV(f_l) \\
  &\leq\lim\limits_{l\to\infty}\left( \frac{1}{2}\|u(f)-u_l^\delta\|_{ L^2(\omega\times(0,T))}^2  + \gamma TV(f) \right) \\
&= \frac{1}{2}\|u(f)-u^\delta\|_{ L^2(\omega\times(0,T))}^2  + \gamma TV(f)=\mathcal{J}(f),
\end{align*}
which implies that $f^*$ is a minimizer of optimization problem $(\mathcal{P})$.
Let us show next the second convergence statement in \eqref{thm-stab-cp-eq2}. We have
\begin{align*}
&\frac{1}{2}\|u(f^*)-u^\delta\|_{ L^2(\omega\times(0,T))}^2 + \gamma\limsup\limits_{l\to \infty} TV(f_l)\\
\leq&  \frac{1}{2}\liminf\limits_{l\to \infty}\|u(f_l)-u_l^\delta\|_{L^2(\omega\times(0,T))}^2  + \gamma \limsup\limits_{l\to \infty} TV(f_l)\\
\leq&\limsup\limits_{l\to \infty} \left(\frac{1}{2}\|u(f_l)-u_l^\delta\|_{L^2(\omega\times(0,T))}^2 + \gamma TV(f_l)\right) \\
\leq&\limsup\limits_{l\to \infty} \left(\frac{1}{2}\|u(f^*)-u_l^\delta\|_{L^2(\omega\times(0,T))}^2 + \gamma TV(f^*)\right) \\
=& \frac{1}{2}\|u(f^*)-u^\delta\|_{L^2(\omega\times(0,T))}^2 + \gamma TV(f^*),
\end{align*}
and hence $\limsup_{l\to \infty} TV(f_l)\leq TV(f^*)$. This together with $TV(f^*)\leq\liminf_{l\to\infty}TV(f_l)$ leads to $TV(f_l)\to TV(f^*)$ as $l\to\infty$.
\hfill$\Box$

\section{Discretization  and convergence}\label{sec:4}

In this section, we propose and analyze a fully discrete method to approximate the nonlinear optimization problem $(\mathcal{P})$.
\subsection{Discretization of the direct problem}
We first propose a full discretization method for the direct problem \eqref{prob-forward}.
The approximation is based on a spatial finite element discretization and a temporal finite difference scheme.
Let $\{\mathcal{T}_h\}_{0<h<1}$ be a family of regular triangulations of the domain $\Omega$, with the mesh parameter $h$ denoting the maximum diameter of the elements.
We introduce the standard finite element spaces  of  piecewise linear functions:
\begin{eqnarray*}
 && U_h:=\left\{v_h\in H^1(\Omega)~|~ v_h|_K\in\mathds{P}_1(K),~\forall K\in\mathcal{T}_h \right\},\\
&& V_h^1:=\left\{v_h\in C(\overline{\Omega})~|~ v_h|_K\in\mathds{P}_1(K),~\forall K\in\mathcal{T}_h \right\},
\end{eqnarray*}
where $\mathds{P}_1(K)$ denotes the space of linear polynomials on $K$.
For the time discretization, we use the popular L1 scheme to discretize the time-fractional derivative \cite{lin2007finite,sun2006fully}.
Let $t_k:=k\tau$, $k=0,1,\cdots,K_\tau$, where $\tau:=T/K_\tau$ is the time step.
The fractional derivative \eqref{det-caputo-1} at $t=t_{k+1}$ is approximated by
\begin{align}
\partial_t^\alpha u(x,t_{k+1}) &= \frac{1}{\Gamma(1-\alpha)} \sum\limits_{j=0}^k \int_{t_j}^{t_{j+1}} (t_{k+1}-s)^{-\alpha}\frac{\partial u(x,s)}{\partial s} \mathrm{d}s  \nonumber\\
& =\frac{1}{\Gamma(1-\alpha)} \sum\limits_{j=0}^k \frac{u(x,t_{j+1})-u(x,t_j)}{\tau} \int_{t_j}^{t_{j+1}} (t_{k+1}-s)^{-\alpha} \mathrm{d}s + r_\tau^{k+1} \nonumber\\
& =\frac{1}{\Gamma(2-\alpha)\tau^\alpha} \sum\limits_{j=0}^k b_j(u(x,t_{k-j+1})-u(x,t_{k-j})) + r_\tau^{k+1},   \nonumber
\end{align}
where $b_j:=(j+1)^{1-\alpha}-j^{1-\alpha}$ for $j=0,1,\cdots,k$, $r_\tau^{k+1}$ is the local truncation error.
The time discretization error was demonstrated to be of order $2-\alpha$ in \cite{lin2007finite,sun2006fully}, assuming the exact solution possesses twice continuous differentiability in time. However, considering the smoothing property of the subdiffusion equation, this regularity condition is overly restrictive and does not hold even for the direct problem with a smooth initial data.
In \cite{jin2016analysis,jin2018numerical}, the error analysis of the L1 scheme was revisited and a convergence rate of $O(\tau)$ for solutions with low regularity was established under the condition
\begin{equation}\label{cond-sm}
   \int_0^{t_k}(t_k-s)^{\alpha-1}\|\mu'(s)f\|_{L^2(\Omega)}\mathrm{d}s<\infty.
\end{equation}

Let $P_h:L^2(\Omega)\to U_h$ denote the standard $L^2$ projection defined by
\begin{equation*}
  (P_h\varphi,v_h)=(\varphi,v_h),\qquad \forall v_h\in U_h.
\end{equation*}
We consider the fully discrete scheme for the direct problem \eqref{prob-forward}: given $u_h^0=0$, find $u_h^{k+1}\in U_h$ ($k=0,1,\cdots,K_\tau-1$), such that
\begin{equation}\label{full-dis-fp}
  (1+\eta)(u_h^{k+1},v_h) + \eta(\nabla u_h^{k+1},\nabla v_h) = \sum\limits_{j=0}^{k-1}(b_j-b_{j+1}) (u_h^{k-j},v_h) + \eta \mu^{k+1}(f_h,v_h)
\end{equation}
holds for all $v_h\in U_h$, where $\eta:=\Gamma(2-\alpha)\tau^\alpha$, $\mu^k=\mu(t_k)$, and $f_h=P_hf$.

\begin{lemma}\label{lem-uncstable-fuly}$($\cite[Lemma 3.3]{jiang2020numerical}$)$
The fully discrete problem \eqref{full-dis-fp} is unconditional stable with respect to the source function $f$ in the sense that for all $h$ and
$\tau>0$, there  exists a constant $c$ independent of $\tau$ and $h$, such that
\begin{equation}\label{disA-bound}
  \|u_h^k\|_{L^2(\Omega)} + \|\nabla u_h^k\|_{L^2(\Omega)}\leq c \|f\|_{L^2(\Omega)}, \qquad \forall k\geq 1.
\end{equation}
\end{lemma}
The condition \eqref{cond-sm} is automatically guaranteed by the assumption that $f\in L^2(\Omega)$ and $\mu\in C^1[0,T]$. Consequently, we derive the following error estimate, which directly follows from \cite[Theorem 3.6]{jin2016two}.
\begin{theorem}\label{thm-fuly-erest}
Assume that $f\in L^2(\Omega)$ and $\mu\in C^1[0,T]$. Let $u$ and $u_h^k$ be the solutions of \eqref{prob-forward} and \eqref{full-dis-fp}, respectively. Then
\begin{eqnarray*}
\|u(t_k)-u_h^k\|_{L^2(\Omega)} \leq c\|\mu\|_{L^\infty(0,T)}\|f\|_{L^2(\Omega)}(t_k^{\alpha-1}\tau  + h^2|\ln h|^2),\qquad \forall k\geq 1.
\end{eqnarray*}
\end{theorem}

\subsection{Finite element approximation of the inverse problem}
The optimization problem $(\mathcal{P})$ can be discreized by
\begin{equation*}
 \min\limits_{f\in F_{ad}^h} \mathcal{J}_{\tau,h}(f),\quad \mathcal{J}_{\tau,h}(f):=\frac{\tau}{2}\sum_{k=0}^{K_\tau}c_k\|u_h^k(f)-u^\delta(\cdot,t_k)\|_{ L^2(\omega)}^2  + \gamma TV(f),\tag{$\mathcal{P}_{\tau,h}$}
\end{equation*}
where $F_{ad}^h:=F_{ad}\cap V_h^1$ is discreized admissible set,
$u_h^{K_{\tau}}(f)$ is the solution of \eqref{full-dis-fp} corresponding to the source term $f$, and $c_k$ represents the coefficients utilized in the composite trapezoidal rule for time integration over the interval $[0,T]$ (i.e., $c_k=1$,  except for  $c_0=c_{K_\tau}=1/2$).

\begin{theorem}\label{thm-discmin-exist}
Let $\gamma>0$ be the fixed regularization parameters.
For every fixed $h,\tau>0$  and $u^\delta\in L^2(\omega)$,
there exists at least a minimizer to the discrete optimization problem $(\mathcal{P}_{\tau,h})$.
\end{theorem}
\noindent \emph{Proof } The assertion is supported by the standing arguments presented in the proof of Theorem 3.1 and \cite[Theorem 3.1]{jiang2020numerical}, which have been omitted for brevity. \hfill$\Box$

The following two convergence conclusions are instrumental in supporting our main conclusion.
\begin{lemma}$($\cite[Lemma 4.6]{hinze2018identifying} and \cite[Lemma 3.2]{hinze2019finite}$)$\label{lem-wfh}
For an fixed $f\in F_{ad}$ an element $\widehat{f_h}\in F_{ad}^h$ exists such that
\begin{equation}\label{lem-Fadh}
  \|\widehat{f_h}-f\|_{L^1(\Omega)}\leq ch|\ln h| \quad\text{and}\quad \lim\limits_{h\to 0}TV(\widehat{f_h})=TV(f).
\end{equation}
\end{lemma}

\begin{lemma}$($\cite[Lemma 4.3]{jiang2020numerical}$)$\label{lem-conres}
Let $\{f_h\}_{h>0}$ be a sequence in $S_h$ which converges weakly to some $f\in L^2(\Omega)$ as $h\to 0$. Let $\{u_h^k(f_h)\}$ and $u(f)$ be the solutions of \eqref{full-dis-fp} and \eqref{prob-forward}, respectively. Then
\begin{equation}\label{eq-convres}
  \lim\limits_{\tau,h\to 0}\frac{\tau}{2}\sum_{k=0}^{K_\tau}c_k\|u_h^k(f_h)-u^\delta(\cdot,t_k)\|_{ L^2(\omega)}^2=\frac{1}{2}\|u(f)-u^\delta\|_{ L^2(\omega\times(0,T))}^2.
\end{equation}
\end{lemma}

The convergence of the finite element approximation ($\mathcal{P}_{\tau,h}$) to the continuous minimization problem ($\mathcal{P}$) is now presented.
\begin{theorem}\label{thm-convfm-ipthm}
Assume that $u^\delta\in L^2(\omega\times(0,T))$. Let $\{h_m\}_{m=1}^\infty$ and $\{\tau_m\}_{m=1}^\infty$ be positive zero sequences as $m\to \infty$,
$f_{m}\in F_{ad}^{h_m}$ is an arbitrary minimizer of the discrete optimization problems $(\mathcal{P}_{\tau,h})$ with $\tau=\tau_m$ and $h=h_m$.
Then an unrelabeled subsequence of $\{f_{m}\}_{m=1}^\infty$ and an element $f^*\in F_{ad}$ exist such that
\begin{eqnarray}
   && \lim\limits_{m\to \infty}\|f_{m}-f^*\|_{L^1(\Omega)}=0 \quad\text{and}\quad \lim\limits_{m\to \infty}TV(f_{m})=TV(f^*). \label{thm-convfm-resipthm}
\end{eqnarray}
Furthermore, $f^*$ is a solution to the continuous problem $(\mathcal{P})$.
\end{theorem}
\noindent \emph{Proof }
Let $f\in F_{ad}$ be arbitrary and $\widehat{f}_{h_m}\in F_{ad}^{h_m}$ be generated from $f$ according to Lemma \ref{lem-wfh}.
The minimizers $f_{m}$ of $(\mathcal{P}_{\tau,h})$ satisfy
\begin{eqnarray*}
   &&  \frac{\tau_m}{2}\sum_{k=0}^{K_{\tau_m}}c_k\|u_{h_m}^{k}(f_{m})-u^\delta(\cdot,t_k)\|_{ L^2(\omega)}^2  + \gamma TV(f_{m}) \\
  && \qquad\qquad\qquad\qquad\leq \frac{\tau_m}{2}\sum_{k=0}^{K_{\tau_m}}c_k\|u_{h_m}^{k}(\widehat{f}_{h_m})-u^\delta(\cdot,t_k)\|_{L^2(\omega)}^2  + \gamma TV(\widehat{f}_{h_m}).
\end{eqnarray*}
By virtue of \eqref{disA-bound} and \eqref{lem-Fadh}, we know that $\{f_{m}\}_{m=1}^\infty$
is  bounded in $BV(\Omega)$-norm. Then it follows from Proposition  \ref{prop-BVfun} that there exists  a (not relabeled) subsequence of
$\{f_{m}\}_{m=1}^\infty$ and an element $f^*\in F_{ad}$ such that
\begin{equation*}
    f_m\to f^*\quad \text{in}\ L^1(\Omega)\  \text{as } m\to \infty,
\end{equation*}
and
\begin{eqnarray}
  TV(f^*)\leq\liminf\limits_{m\to\infty}TV(f_m).\label{thm-convfm-ipthm-eq1}
\end{eqnarray}
Combining this with Lemma \ref{lem-wfh}, Lemma \ref{lem-conres} and the lower semi-continuity of the $L^2$-norm, we arrive at
\begin{align*}
\mathcal{J}(f^*)&=\frac{1}{2}\|u(f^*)-u^\delta\|_{ L^2(\omega\times(0,T))}^2 +  \gamma TV(f^*)  \\
&\leq \lim\limits_{m\to \infty}\frac{\tau_m}{2}\sum_{k=0}^{K_{\tau_m}}c_k\|u_{h_m}^k(f_m)-u^\delta(\cdot,t_k)\|_{ L^2(\omega)}^2 + \gamma\liminf\limits_{m\to \infty} TV(f_m)\\
& \leq \limsup\limits_{m\to\infty} \big(\frac{\tau_m}{2}\sum_{k=0}^{K_{\tau_m}}c_k\|u_{h_m}^k(f_m)-u^\delta(\cdot,t_k)\|_{ L^2(\omega)}^2 + \gamma TV(f_m)\big)\\
& \leq \limsup\limits_{m\to \infty} \big(\frac{\tau_m}{2}\sum_{k=0}^{K_{\tau_m}}c_k\|u_{h_m}^k(\widehat{f}_{h_m})-u^\delta(\cdot,t_k)\|_{ L^2(\omega)}^2 + \gamma TV(\widehat{f}_{h_m})\big)\\
&  =\frac{1}{2}\|u(f)-u^\delta\|_{ L^2(\omega\times(0,T))}^2 +  \gamma TV(f) =\mathcal{J}(f),
\end{align*}
which implies that $f^*$ is a solution to the continuous problem $(\mathcal{P})$.
Finally, the second convergence can be established through a similar discussion as in the proof of Theorem  \ref{thm-stab-cp}.
\hfill$\Box$

The remaining part of this section is dedicated to investigating  the convergence results of a sought source functional as the regularization parameter approaches zero, with an appropriate coupling of noise level and mesh size.
To this end, we consider
\begin{equation*}
  \min\limits_{f\in \mathcal{I}} TV(f), \quad\text{with}\quad \mathcal{I}:=\{f\in F_{ad}~|~ u(f)=u^\dag\ \text{in }\omega\times(0,T) \}. \tag{$\mathcal{IP}$}
\end{equation*}
The existence of the solution to problem ($\mathcal{IP}$) can be derived by the standing arguments of convex analysis.
Moreover, the solution of ($\mathcal{IP}$) can be regarded
as a ``minimal" solution to the unregularized least squares problem
\begin{equation*}
  \min\limits_{f\in F_{ad}} \frac{1}{2}\|u(f)-u^\dag\|_{L^2(\omega)}^2.
\end{equation*}

\begin{theorem}\label{thm-discon}
Let $\{\delta_m\}_{m=1}^\infty$, $\{h_m\}_{m=1}^\infty$, $\{\tau_m\}_{m=1}^\infty$, and $\{\gamma_m\}_{m=1}^\infty$ be positive zero sequences as $m\to \infty$,
and $\{u^{\delta_m}\}_{m=1}^\infty\subset L^2(\omega\times(0,T))$ be a sequence of noisy data satisfying $\|u^{\delta_m}-u^\dag\|_{L^2(\omega\times(0,T))}\leq \delta_m$.
Assume that $f^{\delta_m}$ is an arbitrary minimizer of the discrete optimization problems $(\mathcal{P}_{\tau,h})$
associated to the noisy data $u^{\delta_m}$
with $h=h_m$, $\tau=\tau_m$ and $\gamma=\gamma_m$. Then an unrelabeled subsequence of $\{f^{\delta_m}\}_{m=1}^\infty$ and a total variation minimizing solution $f^\dag$ of the problem $(\mathcal{IP})$ exist such that
\begin{equation*}
  \lim\limits_{m\to\infty}\|f^{\delta_m} -f^\dag\|_{L^1(\Omega)}=0 \quad\text{and}\quad \lim\limits_{m\to\infty}TV(f^{\delta_m})=TV(f^\dag)
\end{equation*}
when
\begin{equation}
   \gamma_m\to 0,\quad \frac{\delta_m}{\sqrt{\gamma_m}}\to 0,
 \quad \frac{h_m|\ln h_m|}{\gamma_m}\to 0, \quad\text{and}\quad \frac{\tau_m}{\sqrt{\gamma_m}}\to 0 \quad\text{as}\quad m\to\infty. \label{thm-discon-cond}
\end{equation}
\end{theorem}
\noindent \emph{Proof }
For any $f\in\mathcal{I}$,
the optimality of $f^{\delta_m}$  implies that
\begin{eqnarray}
  && \frac{\tau_m}{2}\sum_{k=0}^{K_{\tau_m}}c_k\|u_{h_m}^{k}(f^{\delta_m})-u^{\delta_m}(\cdot,t_k)\|_{ L^2(\omega)}^2  + \gamma_m TV(f^{\delta_m}) \nonumber\\
   && \qquad\qquad\qquad \leq \frac{\tau_m}{2}\sum_{k=0}^{K_{\tau_m}}c_k\|u_{h_m}^{k}(\widehat{f}_{h_m})-u^{\delta_m}(\cdot,t_k)\|_{L^2(\omega)}^2  + \gamma_m TV(\widehat{f}_{h_m}).\label{thm-discon-eq1}
\end{eqnarray}
where $\widehat{f}_{h_m}\in F_{ad}^{h_m}$ be generated from $f$ according to Lemma \ref{lem-wfh}.
We bound
\begin{align}
& \frac{\tau_m}{2}\sum_{k=0}^{K_{\tau_m}}c_k\|u_{h_m}^{k}(\widehat{f}_{h_m})-u^{\delta_m}(\cdot,t_k)\|_{L^2(\omega)}^2 \nonumber\\
\leq & \tau_m\sum_{k=0}^{K_{\tau_m}}c_k\|u_{h_m}^{k}(\widehat{f}_{h_m})-u^{\dag}(\cdot,t_k)\|_{L^2(\omega)}^2 + \tau_m\sum_{k=0}^{K_{\tau_m}}c_k\|u^{\dag}(\cdot,t_k)-u^{\delta_m}(\cdot,t_k)\|_{L^2(\omega)}^2 \nonumber\\
\leq & \tau_m\sum_{k=0}^{K_{\tau_m}}c_k\|u_{h_m}^{k}(\widehat{f}_{h_m})-u^{\dag}(\cdot,t_k)\|_{L^2(\omega)}^2 + c\delta_m^2 \label{thm-discon-eq2}
\end{align}
for sufficiently large $m$. Furthermore, applying Lemma \ref{lem-uncstable-fuly}, Theorem \ref{thm-fuly-erest} and the fact that
 \begin{equation*}
  \tau_m\sum_{k=0}^{K_{\tau_m}}c_k\equiv T,\qquad \forall m\in\mathbb{N},
 \end{equation*}
the first term on the right side can be
bounded by
\begin{align*}
& \tau_m\sum_{k=0}^{K_{\tau_m}}c_k\|u_{h_m}^{k}(\widehat{f}_{h_m})-u^{\dag}(\cdot,t_k)\|_{L^2(\omega)}^2 \\
\leq & 2\tau_m\sum_{k=0}^{K_{\tau_m}}c_k\|u_{h_m}^{k}(\widehat{f}_{h_m})-u_{h_m}^{k}(f)\|_{L^2(\omega)}^2 + 2\tau_m\sum_{k=0}^{K_{\tau_m}}c_k\|u_{h_m}^{k}(f)-u(f)(\cdot,t_k)\|_{L^2(\omega)}^2\\
\leq &2\tau_m\sum_{k=0}^{K_{\tau_m}}c_k\left[c\|\widehat{f}_{h_m}-f\|_{L^2(\omega)}^2 +\|u_{h_m}^{k}(f)-u(f)(\cdot,t_k)\|_{L^2(\omega)}^2\right] \\
\leq &c( h_m|\ln h_m|+\tau_m^2 ).   \nonumber
\end{align*}
Combining this with \eqref{thm-discon-eq1} and \eqref{thm-discon-eq2}, we have
\begin{equation*}
  \frac{\tau_m}{2}\sum_{k=0}^{K_{\tau_m}}c_k\|u_{h_m}^{k}(f^{\delta_m})-u^{\delta_m}(\cdot,t_k)\|_{ L^2(\omega)}^2  + \gamma_m TV(f^{\delta_m})\leq c( h_m|\ln h_m|+\tau_m^2+\delta_m^2)+ \gamma_m TV(\widehat{f}_{h_m})
\end{equation*}
for sufficiently large $m$.
Therefore, by using \eqref{lem-Fadh} and \eqref{thm-discon-cond}, we obtain
\begin{eqnarray}
   &&   \lim\limits_{m\to\infty}\frac{\tau_m}{2}\sum_{k=0}^{K_{\tau_m}}c_k\|u_{h_m}^{k}(f^{\delta_m})-u^{\delta_m}(\cdot,t_k)\|_{ L^2(\omega)}^2=0, \label{thm-discon-eq3}\\
   &&  \limsup\limits_{m\to\infty}TV(f^{\delta_m})\leq TV(f). \label{thm-discon-eq4}
\end{eqnarray}
The boundedness of $\{f^{\delta_m}\}_{m=1}^\infty\subset F_{ad}$ can be obtained from \eqref{thm-discon-eq4} as a routine procedure. By virtue of Proposition \ref{prop-BVfun} and Lemma \ref{lem-conres},  a (not relabeled) subsequence of
$\{f^{\delta_m}\}_{m=1}^\infty$  and an element $f^\dag\in F_{ad}$ exist such that
\begin{eqnarray}
   && f^{\delta_m}\to f^\dag\quad \text{in}\ L^1(\Omega)\  \text{as } m\to \infty, \nonumber\\
   &&  TV(f^\dag)\leq\liminf\limits_{m\to\infty}TV(f^{\delta_m}), \label{thm-discon-eq5}
\end{eqnarray}
and
\begin{equation}
    \lim\limits_{m\to \infty}\frac{\tau_m}{2}\sum_{k=0}^{K_{\tau_m}}c_k\|u_{h_m}^k(f^{\delta_m})-u^\dag(\cdot,t_k)\|_{ L^2(\omega)}^2=\frac{1}{2}\|u(f^\dag)-u^\dag\|_{ L^2(\omega\times(0,T))}^2. \label{thm-discon-eq6}
\end{equation}
Then by using \eqref{thm-discon-eq3} and \eqref{thm-discon-eq6}, we have
\begin{eqnarray*}
   &&\frac{1}{2}\|u(f^\dag)-u^\dag\|_{ L^2(\omega\times(0,T))}^2 \\
   &&  \quad \leq \lim\limits_{m\to \infty}\frac{\tau_m}{2}\sum_{k=0}^{K_{\tau_m}}c_k\left(\|u_{h_m}^k(f^{\delta_m}) -u^{\delta_m}(\cdot,t_k)\|_{ L^2(\omega)}^2 + \|u^{\delta_m}(\cdot,t_k)- u^\dag(\cdot,t_k)\|_{ L^2(\omega)}^2\right)=0,
\end{eqnarray*}
which implies that $f^\dag\in\mathcal{I}$. Furthermore,  it follows from \eqref{thm-discon-eq4} and \eqref{thm-discon-eq5} that
\begin{equation}
 TV(f^\dag)\leq\liminf\limits_{m\to\infty}TV(f^{\delta_m})\leq \limsup\limits_{m\to\infty}TV(f^{\delta_m})\leq TV(f). \label{thm-discon-eq7}
\end{equation}
Therefore, $f^\dag$ is a total variation minimizing solution of the problem $(\mathcal{IP})$.
By substituting $f^\dag$ for $f$ in \eqref{thm-discon-eq7}, we can readily derive the second convergence of the conclusion.
\hfill$\Box$

\section{Reconstruction Algorthm}\label{sec:5}

The presence of non-smooth TV terms in the minimization problem poses greater challenges in devising efficient iterative schemes.
Our algorithm is developed based on the TV-dual representation approach \cite{bartels2012total,bartels2015numerical,tian2016linearized,chambolle2011first}, which involves reformulating the minimization problem $(\mathcal{P})$ into a saddle point formulation as follows
\begin{equation}\label{sadd-point-prob}
  \inf\limits_{f\in F_{ad}}\sup\limits_{p\in L^1(\Omega;\mathbb{R}^d)} \Psi(f,p),
\end{equation}
where
\begin{equation*}
  \Psi^\delta(f,p):=\frac{1}{2}\|u(f)-u^\delta\|_{ L^2(\omega\times(0,T))}^2
+ \gamma(\nabla f,p)  -\delta_{\widetilde{\mathcal{B}}_1}(p),
\end{equation*}
$\delta_{\mathcal{B}_1}$ denotes the indicator function of the set $ \widetilde{\mathcal{B}}_1:=\{p\in L^1(\Omega;\mathbb{R}^d)~|~\|p\|_\infty:=\max\limits_{j}\|p_j\|_{L^\infty(\Omega)}\leq 1 \}$, and the inner product
\begin{equation*}
   (\nabla f,p):=\int_{\Omega}\nabla f\cdot p\mathrm{d}x.
\end{equation*}
Using the piecewise constant finite element space
\begin{equation*}
  V_h^0 = \left\{p_h\in L^1(\Omega)~ | ~ p_h|_K=\text{constant},~\forall K\in\mathcal{T}_h \right\},
\end{equation*}
problem $(\mathcal{P})$ is then approximated in the $(F_{ad}^h,V_h^0)$-space by
\begin{equation}\label{dis-sadd-point-prob}
 \inf\limits_{f\in F_{ad}^h}\sup\limits_{p\in (V_h^0)^d} \Psi_{\tau,h}(f,p),
\end{equation}
where
\begin{equation*}
  \Psi_{\tau,h}(f,p):=\frac{\tau}{2}\sum_{k=0}^{K_\tau}c_k\|u_h^k(f)-u^\delta(\cdot,t_k)\|_{ L^2(\omega)}^2  + \gamma(\nabla f,p)- \delta_{\mathcal{B}_1}(p),
\end{equation*}
where $ \mathcal{B}_1:=\{p\in (V_h^0)^d~|~\|p\|_\infty\leq 1 \}$.
The existence of a saddle point in \eqref{sadd-point-prob} and \eqref{dis-sadd-point-prob} can be inferred from
the fact that the objective function is a lower-semicontinuous, proper, convex-concave function, cf., e.g., \cite{rockafellar1997convex} for details.

Before presenting the optimality condition of problem \eqref{dis-sadd-point-prob}, it is imperative to highlight the following two relationships.
\begin{remark}
For each $f\in V_h^1$ the relations
\begin{equation}\label{rel-TV-1}
   TV(f)=\int_{\Omega}|\nabla f|\mathrm{d}x=\sup\limits_{p\in \mathcal{B}_1}
 (\nabla f,p),
\end{equation}
and
\begin{equation}\label{rel-TV-2}
  \partial TV(f)=\left\{ p\in \mathcal{B}_1\subset(V_h^0)^d ~|~ (\nabla f,p)=\int_{\Omega}|\nabla f|\mathrm{d}x\right\}
\end{equation}
hold, where $\partial TV(f)$ is the subgradient of $TV(f)$. For further details, we refer to \cite{bartels2012total} and \cite[Lemma 4.2]{hinze2019finite}.
\end{remark}

Consider the following adjoint problem
\begin{equation}\label{adj-prob}
  \left\{\begin{array}{ll}
               (D_{T^-}^\alpha-\Delta+1 )w(x,t)=\chi_\omega\cdot(u(f)-u^\delta), & (x,t)\in \Omega\times (0,T), \\
               I_{T^-}^{1-\alpha}w(x,t)=0,& (x,t)\in\Omega\times\{T\}, \\
               \partial_\nu w(x,t)=0, & (x,t)\in\partial\Omega\times(0,T),
             \end{array}\right.
\end{equation}
where $\chi_\omega$ denotes the characterization function of $\omega$.
\begin{lemma}$($\cite[Lemma 5.1]{jiang2020numerical}$)$\label{lem-adj}
Let $f,z\in L^2(\Omega)$ and $u(f),w(f)$ be the solutions of systems \eqref{prob-forward} and \eqref{adj-prob} respectively. Then there holds
\begin{equation}\label{rel-adj}
  \int_0^T\int_\omega (u(f)-u^\delta)u(z)\mathrm{d}x\mathrm{d}t= \int_0^T\int_\Omega z\mu w(f)\mathrm{d}x\mathrm{d}t.
\end{equation}
\end{lemma}
\noindent The discrete version of \eqref{rel-adj} is
\begin{equation}\label{dis-rel-adj}
  \tau\sum_{k=0}^{K_\tau}c_k\int_\omega(u_h^k(f)-u^\delta(\cdot,t_k))u_h^k(z)\mathrm{d}x=\tau\sum_{k=0}^{K_\tau}c_k\mu^k(w_h^k(f),z),
\end{equation}
where $u_h^k(f)$ and $w_h^k(f)$ share the same grid. The numerical calculation of equation \eqref{adj-prob} can be performed by setting $\widetilde{w}(x,t)=w(x,T-t)$, resulting in the transformation of  \eqref{adj-prob} into
\begin{equation}\label{w-adj-prob}
  \left\{\begin{array}{ll}
               (D_{0^+}^\alpha-\Delta+1 )\widetilde{w}(x,t)=\chi_\omega\cdot((u(f)-u^\delta)(\cdot,T-t)), & (x,t)\in \Omega\times (0,T), \\
                \lim\limits_{t\to 0^+}I_{0^+}^{1-\alpha}\widetilde{w}(x,t)=0,& x\in\Omega, \\
               \partial_\nu \widetilde{w}(x,t)=0, & (x,t)\in\partial\Omega\times(0,T).
             \end{array}\right.
\end{equation}
By using \eqref{rel-RandC}, the first equation of \eqref{w-adj-prob} can be written as
\begin{equation*}
   (\partial_t^{\alpha}-\Delta+1)\widetilde{w}(x,t)=\chi_\omega\cdot((u(f)-u^\delta)(\cdot,T-t))-\frac{\widetilde{w}(x,0)t^{-\alpha}}{\Gamma(1-\alpha)}.
\end{equation*}
Hence the equation \eqref{adj-prob} can be numerically computed using a discrete method similar to $u_h^k(f)$.  Anyway, we solve the problem \eqref{w-adj-prob} instead of \eqref{adj-prob}, and
then have $w_h^k(f)=\widetilde{w}_h^{K_\tau-k}(f)$ for $k=0,1,\cdots,K_\tau$.

The subsequent lemma presents the first-order optimality condition for problem \eqref{dis-sadd-point-prob}.
\begin{lemma}\label{lem-optimality}$($Optimality$)$
The function $f\in F_{ad}^h$ is a solution of \eqref{dis-sadd-point-prob} if and only if there exists $p\in \partial TV(f)$ such that
\begin{eqnarray}
   && \tau\sum_{k=0}^{K_\tau}c_k\mu^k(w_h^k(f),g-f)+\gamma (p,\nabla (g-f)) \geq 0, \qquad \forall g\in F_{ad}^h, \label{lem-opti-cond1}\\
   && (\nabla f,q-p)\leq 0, \qquad  \forall q\in \mathcal{B}_1,  \label{lem-opti-cond2}
\end{eqnarray}
\end{lemma}
\noindent \emph{Proof }
The optimality condition is derived from the standard argument (see, e.g., \cite[Lemma 2.21]{troltzsch2010optimal} and \cite[Lemma 10.3]{bartels2015numerical}), we have
\begin{equation}\label{lem-optimality-eq1}
 \tau\sum_{k=0}^{K_\tau}c_k\int_\omega(u_h^k(f)-u^\delta(\cdot,t_k))u_h^k(g-f)\mathrm{d}x+\gamma\int_{\Omega}\nabla(g-f)\cdot p\mathrm{d} x\geq 0,\quad \forall g\in F_{ad}^h,
\end{equation}
and
\begin{equation}\label{lem-optimality-eq2}
  0\in \partial_\rho \Psi_{\tau,h}(f,p).
\end{equation}
According to \eqref{dis-rel-adj}, we obtain \eqref{lem-opti-cond1} from \eqref{lem-optimality-eq1}.

For the second part of the conclusion, \eqref{lem-optimality-eq2} implies that there exists $\xi\in\partial\delta_{\mathcal{B}_1}(p)$ such that $\xi=\nabla f$.
This together with
\begin{equation*}
   \partial\delta_{\mathcal{B}_1}(p)=\left\{ \xi\in((V_h^0)^d)^* ~|~ \langle\xi,q-p\rangle \leq 0 \quad\text{for~all }q\in{\mathcal{B}_1} \right\}
\end{equation*}
yields \eqref{lem-opti-cond2}. Finally, according to \eqref{rel-TV-1} and \eqref{lem-opti-cond2},we arrive at
\begin{equation*}
 (\nabla f,p)=\int_{\Omega}|\nabla f|\mathrm{d}x.
\end{equation*}
Hence it follows from \eqref{rel-TV-2} that $p\in\partial TV(f)$.
\hfill$\Box$

\begin{remark}\label{rem-optimality}
We can rewrite the system \eqref{lem-opti-cond1}-\eqref{lem-opti-cond2} as the following variational inequality in a compact form:
find $\nu\in F_{ad}^h\times \partial TV(f)$, such that
\begin{eqnarray}
   && (F(\nu),\kappa-\nu)\geq 0, \qquad \forall \kappa\in F_{ad}^h\times \mathcal{B}_1, \label{rem-optimality-comp}
\end{eqnarray}
where
\begin{equation*}
  \nu:=\left(\begin{array}{c}
               f \\
               p
             \end{array}\right),\quad \kappa:=\left(\begin{array}{c}
               g \\
               q
             \end{array}\right),\quad F(\nu):=\left(\begin{array}{c}
                      \tau\sum\limits_{k=0}^{K_\tau}c_k\mu^kw_h^k(f)-\gamma\mathrm{div}p\\
                      -\nabla f
                    \end{array}\right),
\end{equation*}
and $-\mathrm{div}$ is the adjoint operators of  $\nabla$ defined by
\begin{equation*}
  (-\mathrm{div}q,g)=(q,\nabla g),\qquad \forall q\in(V_h^0)^d,\quad \forall g\in V_h^1.
\end{equation*}
\end{remark}

The saddle point problem \eqref{dis-sadd-point-prob} can be efficiently solved using the primal-dual algorithm.
The primal-dual algorithm is a popular method for solving optimization problems, aiming to iteratively approach the optimal solution by alternating between optimizing the primal and dual variables.
In terms of efficiency, the primal-dual algorithm offers several  advantages over other optimization methods. It can handle large-scale problems efficiently due to its iterative nature and ability to decompose the problem into smaller subproblems. Moreover, it effectively addresses problems with complex structures and nonlinearities due to its adaptability to different types of constraints.

In this article, we employ a linearized primal-dual algorithm, inspired by the concept presented in \cite{tian2016linearized}, to tackle the saddle point problem \eqref{dis-sadd-point-prob}.

\begin{algorithm} $($Linearized primal-dual algorithm$)$
\label{alg-acc-pd}

\noindent -- Input: Let parameters $\gamma$, $\varsigma$, $\theta>0$ such that
\begin{equation}
  \left(\frac{1}{\varsigma}-c^2\right)\frac{\theta}{\varsigma} >\gamma^2\|\nabla\|^2,\label{cond1-AG-conv}
\end{equation}
where $c$ is the norm of the operator $\mathcal{A}:f\mapsto u(f)$ from $L^2(\Omega)$ to $L^2(\omega\times (0,T))$. Choose an initial guess
$(f^0,p^0)\in F_{ad}^h\times (V_h^0)^d$.

\noindent -- For $n=0,1,2\cdots,$ do

Update the new iteration $(f^{n+1},p^{n+1})$ via solving
\begin{align}
\widetilde{f}^{n+1}:&=\mathop{\arg\min}\limits_{f\in F_{ad}^h}\left\{ \tau\sum_{k=0}^{K_\tau}c_k\mu^k(w_h^k(f^n),f) +\gamma (\nabla f,p^{n}) + \frac{1}{2\varsigma}\|f-f^n\|_{L^2(\Omega)}^2\right\}, \label{iter-pd-c}\\
    f^{n+1}:&=2\widetilde{f}^{n+1}-f^{n},  \label{iter-pd-a} \\
   p^{n+1}:&=\mathop{\arg\max}\limits_{\rho\in (V_h^0)^d}\left\{ \gamma (\nabla f^{n+1},p^n) -\delta_{\mathcal{B}_1}(\rho^n) -\frac{\theta}{2\varsigma}\|p-p^n\|_{L^2(\Omega)}^2  \right\}. \label{iter-pd-b}
\end{align}

\noindent \ \ end

\noindent -- Output: An approximation solution $f^n$ if $n\leq N_{\max}$ for some $N_{\max}\in\mathds{N}$.

\end{algorithm}

\begin{remark}
$(1)$ The notation
\begin{equation*}
  \|\nabla\|:=\sup_{0\neq f\in V_h^1} \frac{\|\nabla f\|_{L^2(\Omega)}}{\|f\|_{L^2(\Omega)}},
\end{equation*}
and the well known inverse inequality implies that $\|\nabla\|\leq ch^{-1}$.

$(2)$ The optimality of $\widetilde{f}^{n+1}$  in the step \eqref{iter-pd-c} yields
\begin{equation*}
   \tau\sum_{k=0}^{K_\tau}c_k\mu^k(w_h^k(f^n),g-\widetilde{f}^{n+1})+\gamma (p^n,\nabla (g-\widetilde{f}^{n+1}))+\frac{1}{\varsigma}(\widetilde{f}^{n+1}-f^n,g-\widetilde{f}^{n+1}) \geq 0
\end{equation*}
for all $g\in F_{ad}^h$, which can be rewritten as
\begin{equation*}
  (f^n+\varsigma(\gamma\mathrm{div}p^n-\tau\sum_{k=0}^{K_\tau}c_k\mu^k(w_h^k(f^n)))-\widetilde{f}^{n+1} ,g-\widetilde{f}^{n+1}) \leq 0
\end{equation*}
for all $g\in F_{ad}^h$. Hence we have
\begin{equation}\label{exfor-PD-f}
  \widetilde{f}^{n+1}=\mathbb{P}_{F_{ad}^h}(f^n+\varsigma(\gamma\mathrm{div}p^n-\tau\sum_{k=0}^{K_\tau}c_k\mu^k(w_h^k(f^n))),
\end{equation}
where $\mathbb{P}_{F_{ad}^h}:L^2(\Omega)\to F_{ad}^h$ denotes the $L^2$-projection on the set $F_{ad}^h$ characterized by
\begin{equation*}
 (\phi-\mathbb{P}_{F_{ad}^h}\phi ,g-\mathbb{P}_{F_{ad}^h}\phi) \leq 0
\end{equation*}
for all $g\in F_{ad}^h$.

$(3)$ The solution $\rho^{n+1}$ to the subproblem \eqref{iter-pd-b} is explicitly given by
\begin{eqnarray}
   && p^{n+1}=\frac{p^n+\frac{\gamma\varsigma}{\theta}\nabla\widetilde{f}^{n+1} }{\max\big\{ 1,|p^n+\frac{\gamma\varsigma}{\theta}\nabla\widetilde{f}^{n+1}| \big\}},  \label{exfor-PD-rho}
\end{eqnarray}
which can be computed element-wise $($refer to \cite{bartels2012total,bartels2015numerical,tian2016linearized}$)$.

$(4)$ The $O(1/n)$ convergence of the Algorithm \ref{alg-acc-pd} can be derived from \cite{tian2016linearized}, where $n$  denotes the iteration counter.
\end{remark}

\section{Numerical results}\label{sec:6}
In this section, we present some numerical results to show the efficiency of model ($\mathcal{P}_{\tau,h}$) and Algorithm \ref{alg-acc-pd}.
All these experiments are implemented in Matlab R2021b and run on a personal computer with 2.80 GHz CPU processor.

The exact data $u^\dag=u(f^*)|_{\omega}$ is obtained in our numerical tests by solving the direct problem \eqref{prob-forward} with the true source function $f^*$. The noisy data $u^\delta$ is then generated by adding a random uniform perturbation on $[-1,1]$, i.e.,
\begin{equation*}
  u^\delta(x,t)=u^\dag(x,t)+\delta_{\mathrm{rel}} {r\over\|r\|_{L^2(\omega\times(0,T))}}u^\dag(x,t),\quad (x,t)\in\omega\times(0,t]
\end{equation*}
where $r(\cdot,t)=2\text{rand}(\text{size}(u^\dag(\cdot,t)))-1$ and $\delta_{\mathrm{rel}}>0$ is the relative noise level.
Then the corresponding noise level in \eqref{noisy-data} is calculated by $\delta:=\delta_{\mathrm{rel}}\|u^\dag\|_{L^2(\omega\times(0,T))}$.
The proposed algorithm is validated by investigating the behavior of the reconstructed source error and its corresponding residue, as defined by
\begin{equation*}
   e_r(f^n,f^*):=\frac{\|f^n-f^*\|_{L^2(\Omega)} }{\|f^*\|_{L^2(\Omega)}} \quad\text{and}\quad \mathrm{res}(f^n,u^\delta):=\tau\sum_{k=0}^{K_\tau}c_k\|u_h^k(f^n)-u^\delta(\cdot,t_k)\|_{ L^2(\omega)}^2.
\end{equation*}
Inspired by \cite{tian2016linearized} and the concept of iterated Tikhonov regularization \cite{hanke1998nonstationary}, iterative algorithm \ref{alg-acc-pd} terminates either when the number of iterations n reaches the maximum iteration count $N_{\max}=3000$ or when $n$ is the minimum value such that
\begin{equation*}
  \frac{\|f^{n}-f^{n-1}\|_{L^2(\Omega)}}{\|f^{n}\|_{L^2(\Omega)} }\leq 10^{-4} \quad\text{or}\quad    \mathrm{res}(f^n,u^\delta)\leq 1.01\delta,
\end{equation*}
where the second stopping criterion is classical Morozov discrepancy principle (see, e.g., \cite[Chapter 4]{engl1996regularization}).

We first consider a one-dimensional problem with region $\overline{\Omega}\times[0,T]=[0,1]^2$ and test Algorithm \ref{alg-acc-pd} with exact source function $f^*$ as follows:
\begin{example}\label{examp-1}
$f^*(x)=0.5\chi_{[1/4,3/4]}$.
\end{example}
\noindent In our numerical test, the equation is discretized by dividing $[0,1]^2$ into $50\times 50$ equidistant meshes.
The selected observation areas are as follows: (1) $\omega=[2/50,25/50]$; (2) $\omega=[28/50,48/50]$; (3) $\omega=[2/50,25/50]\cup [28/50,48/50]$.
The control variable parameters are set as $\underline{f}=0$ and $\overline{f}=1$.
For known  temporal distribution, we consider $\mu(t)\equiv 1$ and $\mu(t)=\cos(2\pi t)$ respectively.
Parameters $\zeta$ and $\theta$ must be chosen carefully so that condition \eqref{cond1-AG-conv} is satisfied.
After performing preliminary calculations and averaging,
we have $\|u(f^*)\|_{ L^2(\Omega\times(0,T))}/\|f^*\|_{L^2(\Omega)}\approx 0.0623$ for $\alpha=0.3$ and $\|u(f^*)\|_{ L^2(\Omega\times(0,T))}/\|f^*\|_{L^2(\Omega)}\approx 0.0606$ for $\alpha=0.8$,
thus we set accordingly $\varsigma=100$ and $\theta=10^{-1}$ in all the one-dimensional tests.
In fact, the flexibility of choosing hyperparameters $\varsigma$ and $\theta$ is much greater due to the typically smaller scale of the regularization parameter $\gamma$.
As an initial approximation we choose $f^0(x)\equiv 0.25$, and $p^0(x)\equiv 0.5$.

\begin{figure}[h!]
 \centerline{\scalebox{0.38}{\includegraphics{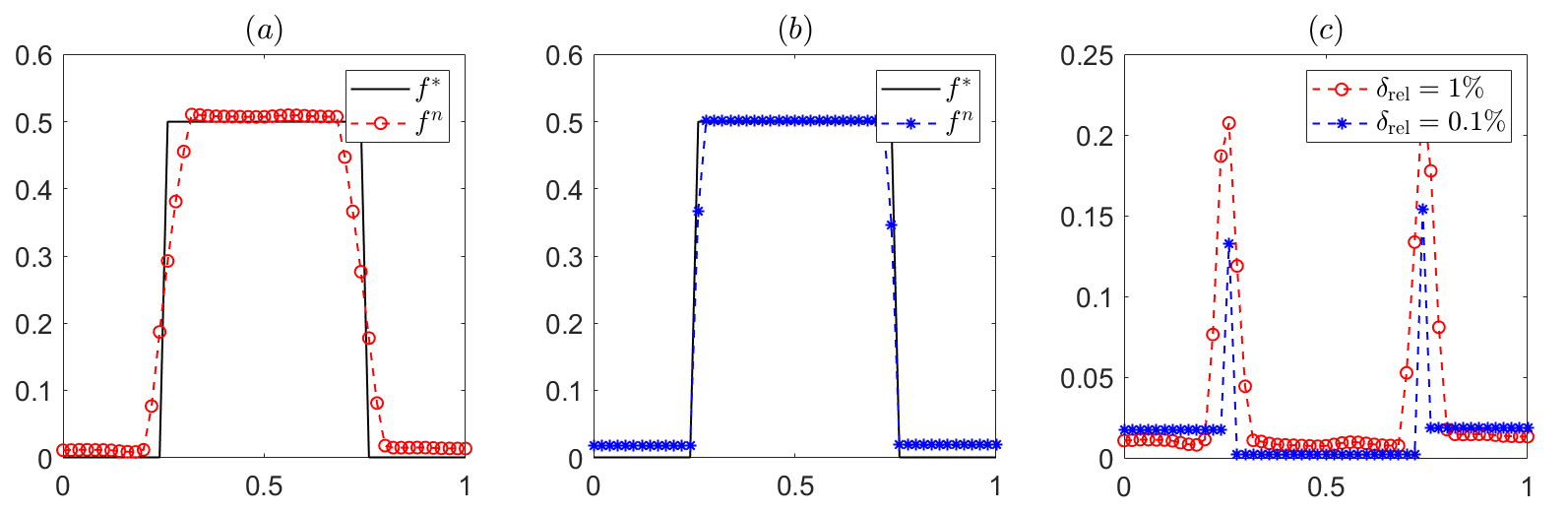}}}
  \centering\caption{\small The computed source functions for Example \ref{examp-1} with $\alpha=0.3$, $\mu(t)=\cos(2\pi t)$ and $\omega=[2/50,25/50]$. $(a)$: $\delta_{\mathrm{rel}}=1\%$; $(b)$: $\delta_{\mathrm{rel}}=0.1\%$; $(c)$: the absolute error $|f^n-f^*|$.} \label{Fig-1D-1}
\end{figure}

\begin{figure}[h!]
 \centerline{\scalebox{0.38}{\includegraphics{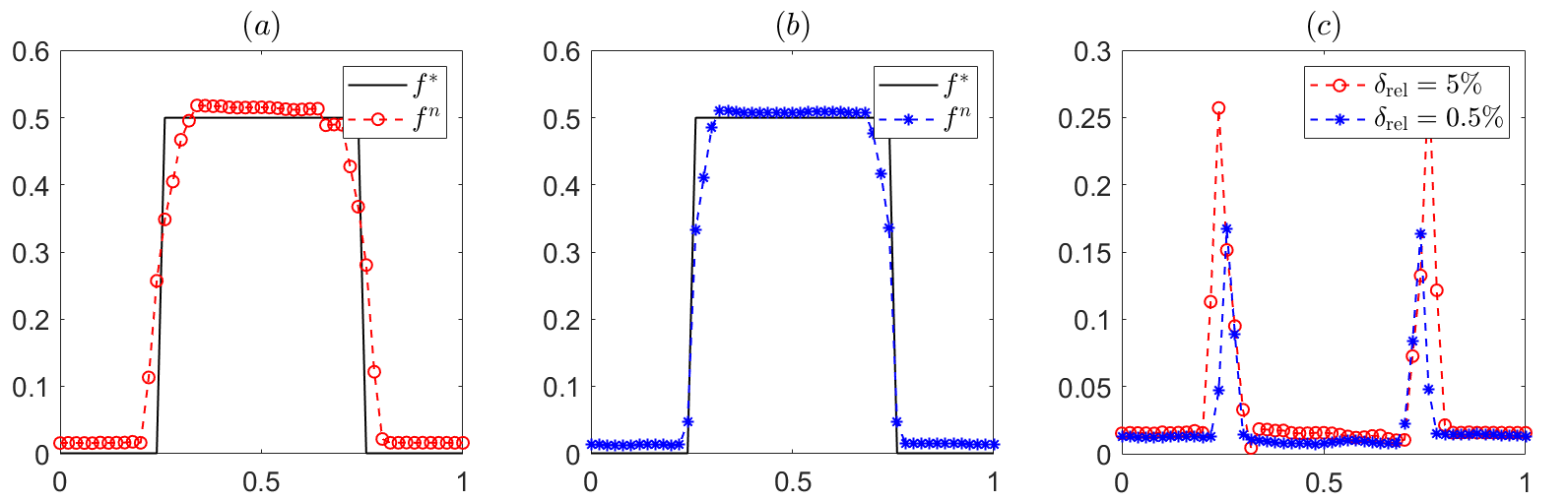}}}
  \centering\caption{\small The computed source functions for Example \ref{examp-1} with $\alpha=0.8$, $\mu(t)\equiv 1$ and $\omega=[28/50,48/50]$. $(a)$: $\delta_{\mathrm{rel}}=5\%$; $(b)$: $\delta_{\mathrm{rel}}=0.5\%$; $(c)$: the absolute error $|f^n-f^*|$.} \label{Fig-1D-2}
\end{figure}

The obtained numerical results are listed in Table \ref{Tab1} and \ref{Tab2}, including the impact of noise level, regularization parameters, and time-fractional order on the reconstruction quality.  In general, the regularized solution becomes more accurate and the number of iterations increases as the noise level of data decreases.
Figure \ref{Fig-1D-1} and \ref{Fig-1D-2} display some of results for the regularized solutions and absolute errors, which show that our algorithm is effective.

\begin{table}[h!]\small
\setlength{\belowcaptionskip}{0.2cm}
\centering
\caption{Numerical results of the reconstruction for Example \ref{examp-1} with $\mu(t)=\cos(2\pi t)$}\label{Tab1}\setlength{\belowcaptionskip}{-0.9cm}
\begin{tabular}{ccccccccc}
\hline
   &  &  &  \multicolumn{3}{c}{$\alpha=0.3$ } &  \multicolumn{3}{c}{$\alpha=0.8$ }  \\
   \cmidrule(lr){4-6} \cmidrule(lr){7-9}
$\delta_{\mathrm{rel}}$  & $\omega$  &   $\gamma$  & $n$  & $e_r(f^n,f^*)$  &  $\mathrm{res}(f^n,g^\delta)$   & $n$  & $e_r(f^n,f^*)$  &  $\mathrm{res}(f^n,g^\delta)$ \\\hline
$0.1\%$  & (1)    & $10^{-9}$   &  $862$  &  $0.0382$   & $1.2067e$--$05$    & $843$   &  $0.0441$   & $1.0367e$--$05$   \\
   & (2) &   $10^{-9}$  & $668$  &  $0.0489$   &   $1.3398e$--$05$  &  $712$  &  $0.0388$   &  $1.1265e$--$05$  \\
   & (3) &   $10^{-9}$  & $565$   &  $0.0567$   & $1.9454e$--$05$    & $572$  & $0.0503$   & $1.7408e$--$05$   \\ \hline
$0.5\%$  &  (1)   & $10^{-9}$   & $372$  &  $0.0819$  &  $5.7258e$--$05$  & $320$  &  $0.0914$  &  $5.7216e$--$05$ \\
   & (2) &   $10^{-9}$  & $312$  &  $0.0923$  &  $6.0147e$--$05$  & $286$  &  $0.1023$  &  $5.3218e$--$05$ \\
   & (3) &   $10^{-9}$  & $322$  &  $0.0823$  &  $9.2306e$--$05$  & $298$  &  $0.1054$  &  $8.4235e$--$05$ \\ \hline
$1\%$  &  (1)   & $10^{-8}$   & $257$  &  $0.1119$  &  $1.0121e$--$04$  & $241$  &  $0.1391$  &  $1.0121e$--$04$ \\
   & (2) &   $10^{-8}$  & $195$  &  $0.1343$  &  $1.2017e$--$04$  & $219$  &  $0.1322$  &  $1.0912e$--$04$ \\
   & (3)  &   $10^{-8}$  & $295$  &  $0.1298$  &  $1.9043e$--$04$  & $281$  &  $0.1308$  &  $1.6824e$--$04$ \\ \hline
   $5\%$  &  (1)   & $10^{-8}$   & $211$  &  $0.1760$  &  $5.0134e$--$04$  & $208$  &  $0.1703$  &  $5.0112e$--$04$ \\
   & (2) &   $10^{-8}$  & $176$  &  $0.1781$  &  $6.8530e$--$04$  & $183$  &  $0.1760$  &  $6.2303e$--$04$ \\
   & (3)  &   $10^{-8}$  & $237$  &  $0.1795$  &  $8.0012e$--$04$  & $239$  &  $0.1798$  &  $5.6451e$--$04$ \\ \hline
\end{tabular}
\end{table}

\begin{table}[h!]\small
\setlength{\belowcaptionskip}{0.2cm}
\centering
\caption{Numerical results of the reconstruction for Example \ref{examp-1} with $\mu(t)\equiv 1$}\label{Tab2}\setlength{\belowcaptionskip}{-0.9cm}
\begin{tabular}{ccccccccc}
\hline
   &  &  &  \multicolumn{3}{c}{$\alpha=0.3$ } &  \multicolumn{3}{c}{$\alpha=0.8$ }  \\
   \cmidrule(lr){4-6} \cmidrule(lr){7-9}
$\delta_{\mathrm{rel}}$  & $\omega$  &   $\gamma$  & $n$  & $e_r(f^n,f^*)$  &  $\mathrm{res}(f^n,g^\delta)$   & $n$  & $e_r(f^n,f^*)$  &  $\mathrm{res}(f^n,g^\delta)$ \\\hline
$0.1\%$  & (1)    & $10^{-9}$   &  $831$  &  $0.0361$   & $1.1872e$--$05$    & $822$   &  $0.0423$   & $1.0987e$--$05$   \\
   & (2) &   $10^{-9}$  & $656$  &  $0.0462$   &   $1.3007e$--$05$  &  $691$  &  $0.0386$   &  $1.1102e$--$05$  \\
   & (3) &   $10^{-9}$  & $530$   &  $0.0503$   & $1.8982e$--$05$    & $561$  & $0.0481$   & $1.6838e$--$05$   \\ \hline
$0.5\%$  &  (1)   & $10^{-9}$   & $351$  &  $0.0823$  &  $6.0022e$--$05$  & $314$  &  $0.0885$  &  $5.9242e$--$05$ \\
   & (2) &   $10^{-9}$  & $308$  &  $0.0911$  &  $6.0002e$--$05$  & $282$  &  $0.0988$  &  $5.0890e$--$05$ \\
   & (3) &   $10^{-9}$  & $297$  &  $0.0807$  &  $8.7217e$--$05$  & $289$  &  $0.0965$  &  $5.7821e$--$05$ \\ \hline
$1\%$  &  (1)   & $10^{-8}$   & $241$  &  $0.1028$  &  $1.0007e$--$04$  & $234$  &  $0.1198$  &  $0.9972e$--$04$ \\
   & (2) &   $10^{-8}$  & $198$  &  $0.1245$  &  $1.1882e$--$04$  & $208$  &  $0.1201$  &  $1.0322e$--$04$ \\
   & (3)  &   $10^{-8}$  & $281$  &  $0.1176$  &  $1.3036e$--$04$  & $276$  &  $0.1268$  &  $1.4214e$--$04$ \\ \hline
   $5\%$  &  (1)   & $10^{-8}$   & $209$  &  $0.1709$  &  $5.0002e$--$04$  & $201$  &  $0.1693$  &  $4.8977e$--$04$ \\
   & (2) &   $10^{-8}$  & $185$  &  $0.1766$  &  $5.7934e$--$04$  & $193$  &  $0.1689$  &  $5.9972e$--$04$ \\
   & (3)  &   $10^{-8}$  & $245$  &  $0.1690$  &  $5.1025e$--$04$  & $220$  &  $0.1688$  &  $4.0871e$--$04$ \\ \hline
\end{tabular}
\end{table}

\begin{figure}[h!]
 \centerline{\scalebox{0.38}{\includegraphics{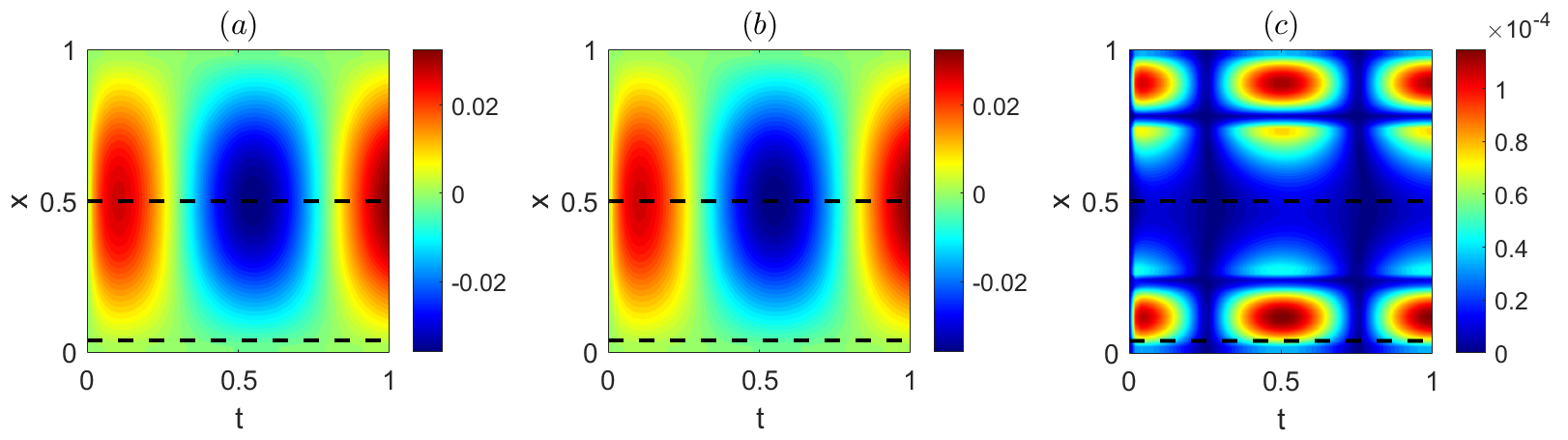}}}\vspace{5pt}
 \centerline{\scalebox{0.38}{\includegraphics{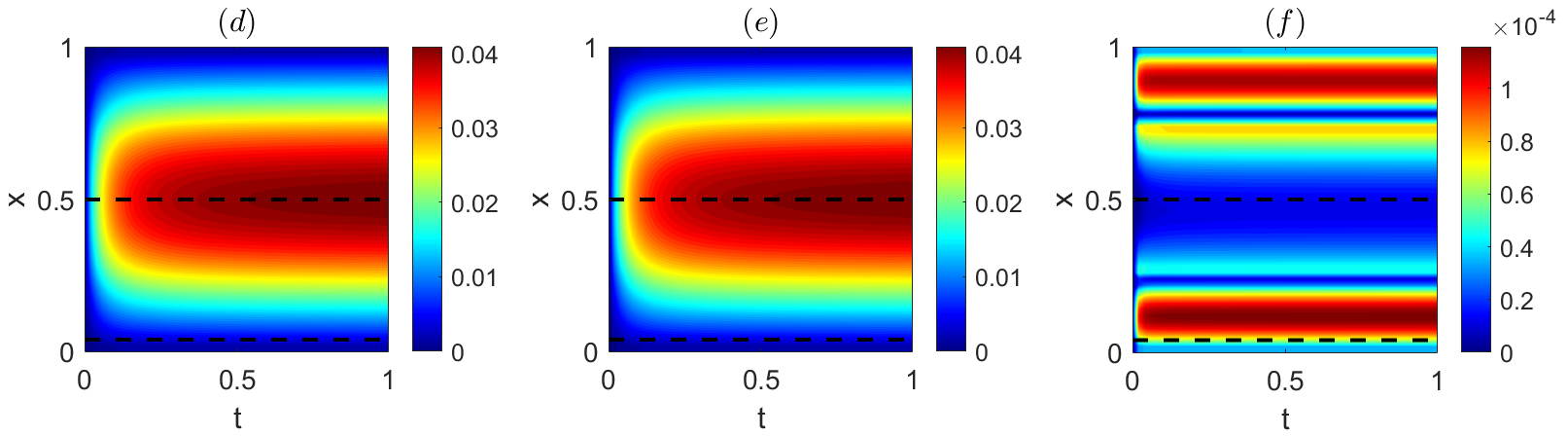}}}
  \centering\caption{\small Approximate solution of equation \eqref{prob-forward} with $\alpha=0.8$ and $\mu(t)=\cos(2\pi t)$(first row)/ $\mu(t)\equiv 1$(second row). The observation regions are located  between the black dashed lines, i.e., $\omega=[2/50,25/50]$. The relative noise level $\delta_{\mathrm{rel}}=0.1\%$. $(a)$ and $(d)$: finite element solution $u(f^*)(x,t)$; $(b)$ and $(e)$: approximation solution $u(f^n)(x,t)$; $(c)$ and $(f)$:  the absolute error $|u(f^n)(x,t)-u(f^*)(x,t)|$.
  } \label{Fig-1D-3}
\end{figure}

We also present the solution of equation \eqref{prob-forward} using the approximate source term we obtained, and compare it with the finite element solution under the exact source, as illustrated in Figure \ref{Fig-1D-3}.
One of the main difference between inverse problem and the direct problem is the way they solve the problem.
In direct problems, we initially establish a well-defined model and subsequently employ various solvers such as the finite element method to solve the corresponding equation.
Conversely, inverse problems enable us to reconstruct the model from data without complete knowledge of its characteristics and provide a solution for the equation.
Therefore, in practical scenarios, the inverse problem model is typically  more flexible than the direct problem model.

In order to provide a more comprehensive evaluation of the impact of TV regularization, we incorporate the classical Tikhonov model as a benchmark for comparison, that is
\begin{equation*}
 \min\limits_{f\in F_{ad}} \mathcal{T}(f),\quad \mathcal{T}(f):=\frac{1}{2}\|u(f)-u^\delta\|_{ L^2(\omega\times(0,T))}^2  + \beta\|f\|_{L^2(\Omega)},\tag{$\mathcal{T}$}
\end{equation*}
where $\beta>0$ is the regularization parameter. Then the discrete version is
\begin{equation*}
 \min\limits_{f\in F_{ad}^h} \mathcal{T}_{\tau,h}(f),\quad \mathcal{T}_{\tau,h}(f):=\frac{\tau}{2}\sum_{k=0}^{K_\tau}c_k\|u_h^k(f)-u^\delta(\cdot,t_k)\|_{ L^2(\omega)}^2  + \beta\|f\|_{L^2(\Omega)}.\tag{$\mathcal{T}_{\tau,h}$}
\end{equation*}
The problem $\mathcal{T}_{\tau,h}$ is solved using an iterative method similar to algorithm \ref{alg-acc-pd}, where the sought term is updated as
\begin{equation*}
  f^{n+1}:=\mathop{\arg\min}\limits_{f\in F_{ad}^h}\left\{ \tau\sum_{k=0}^{K_\tau}c_k\mu^k(w_h^k(f^n),f) +\beta\|f\|_{L^2(\Omega)} + \frac{1}{2\vartheta}\|f-f^n\|_{L^2(\Omega)}^2\right\},
\end{equation*}
The optimality of $f^{n+1}$ in the above step yields
\begin{equation*}
  (f^n-\vartheta(\tau\sum_{k=0}^{K_\tau}c_k\mu^k(w_h^k(f^n)+2\beta f^n)-f^{n+1},g-f^{n+1}) \leq 0
\end{equation*}
for all $g\in F_{ad}^h$. Hence we have
\begin{equation}\label{iter-Tk}
  f^{n+1}=\mathbb{P}_{F_{ad}^h}(f^n-\vartheta(\tau\sum_{k=0}^{K_\tau}c_k\mu^k(w_h^k(f^n)+2\beta f^n)).
\end{equation}
The iterative scheme \eqref{iter-Tk} essentially belongs to a class of projected gradient descent methods, where the stepsize $\vartheta>0$ needs to be carefully chosen. In this test, we set $\alpha=0.5$ and $\mu=\sin(5\pi t)$ for Example \ref{examp-1}. The observation area is $\omega=[15/50,35/50]$ and is subject to noise at a relative error level  $\delta_{\mathrm{rel}}=0.5\%$.
The stepsize $\vartheta$ is set to $0.1$, while maintaining the previous discretization and optimization parameters.

\begin{figure}[h!]
 \centerline{\scalebox{0.38}{\includegraphics{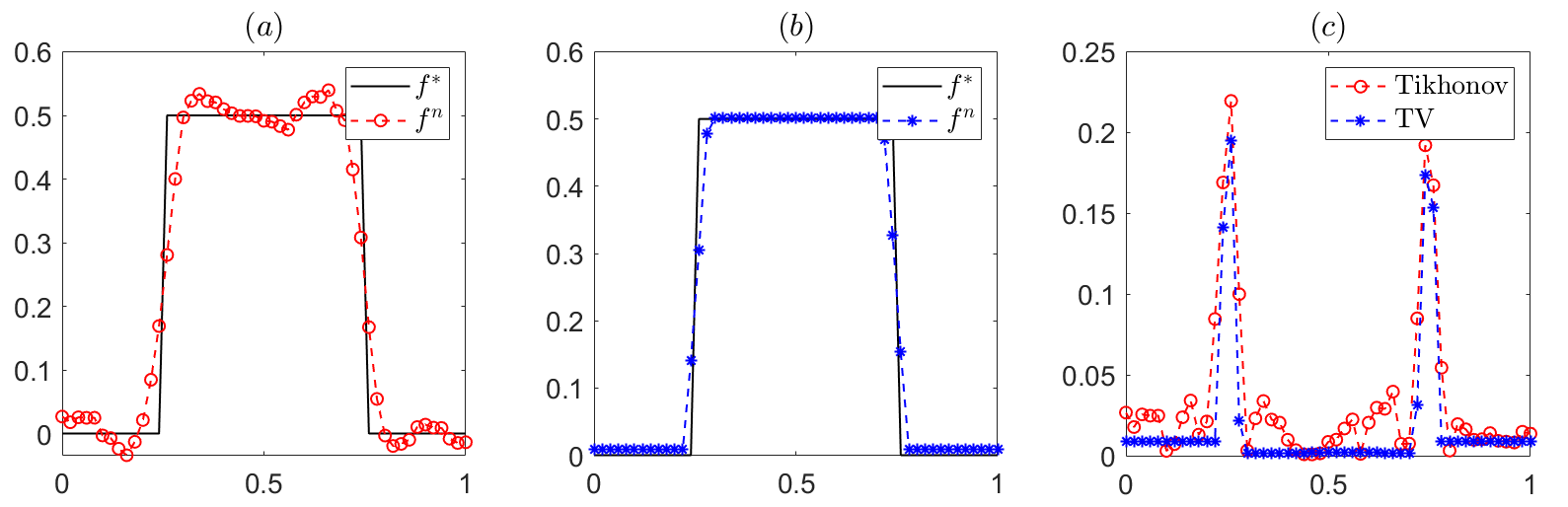}}}
  \centering\caption{\small The computed source functions for Example \ref{examp-1} with $\alpha=0.5$, $\mu(t)=\sin(5\pi t)$, $\omega=[15/50,35/50]$ and $\delta_{\mathrm{rel}}=0.5\%$. $(a)$: Tikhonov solution; $(b)$: TV solution; $(c)$: the absolute error $|f^n-f^*|$.} \label{Fig-1D-4}
\end{figure}

\begin{figure}[h!]
 \centerline{\scalebox{0.38}{\includegraphics{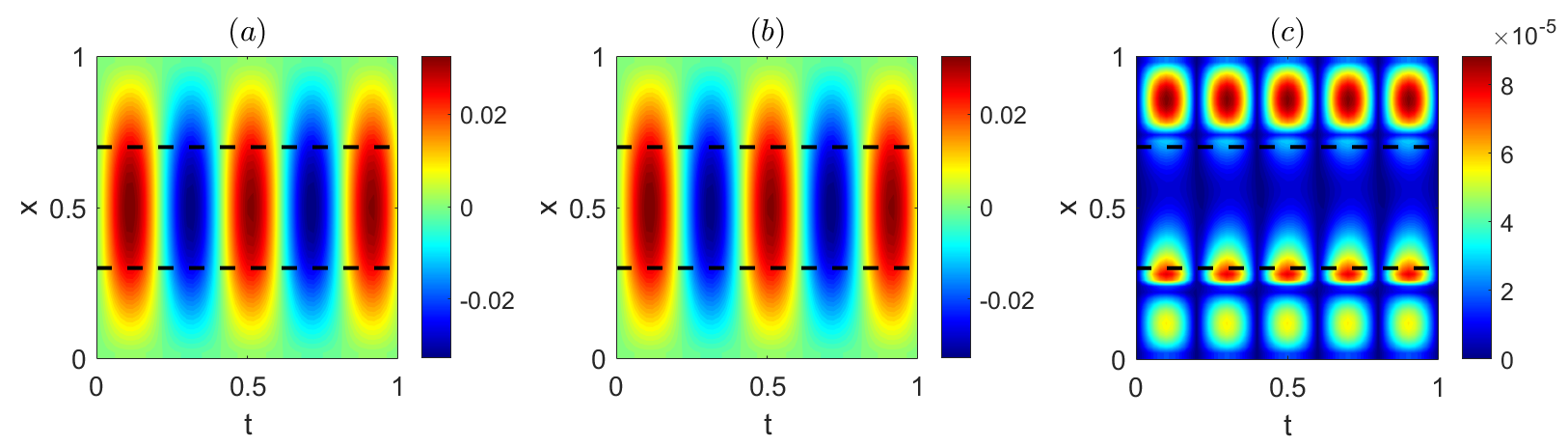}}}\vspace{5pt}
 \centerline{\scalebox{0.38}{\includegraphics{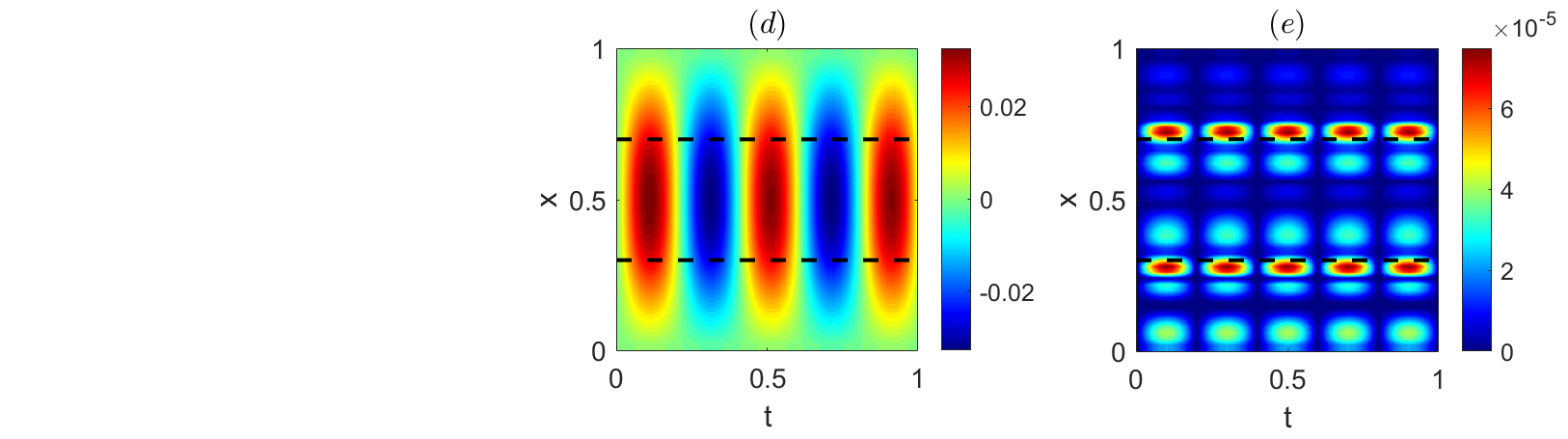}}}
  \centering\caption{\small Approximate solution of equation \eqref{prob-forward} with $\alpha=0.5$ and $\mu(t)=\sin(5\pi t)$. The observation regions are located between the black dashed lines, i.e., $\omega=[15/50,35/50]$. The relative noise level $\delta_{\mathrm{rel}}=0.1\%$.
  $(a)$: finite element solution $u(f^*)(x,t)$; $(b)$: Tikhonov solution $u(f^n)(x,t)$ $(c)$:  the
absolute error of Tikhonov solution; $(d)$: TV solution; $(c)$:  the
absolute error of TV solution.
  } \label{Fig-1D-5}
\end{figure}

The test results are depicted in Figures \ref{Fig-1D-4} and \ref{Fig-1D-5}.
It is not difficult to see that the TV regularization term has the effect of stabilizing the place where
the gradient changes greatly. In other words, the proposed method effectively eliminates high-frequency oscillations, thereby enhancing the accuracy of solution recovery near the discontinuous points.

\begin{figure}[h!]
 \centerline{\scalebox{0.36}{\includegraphics{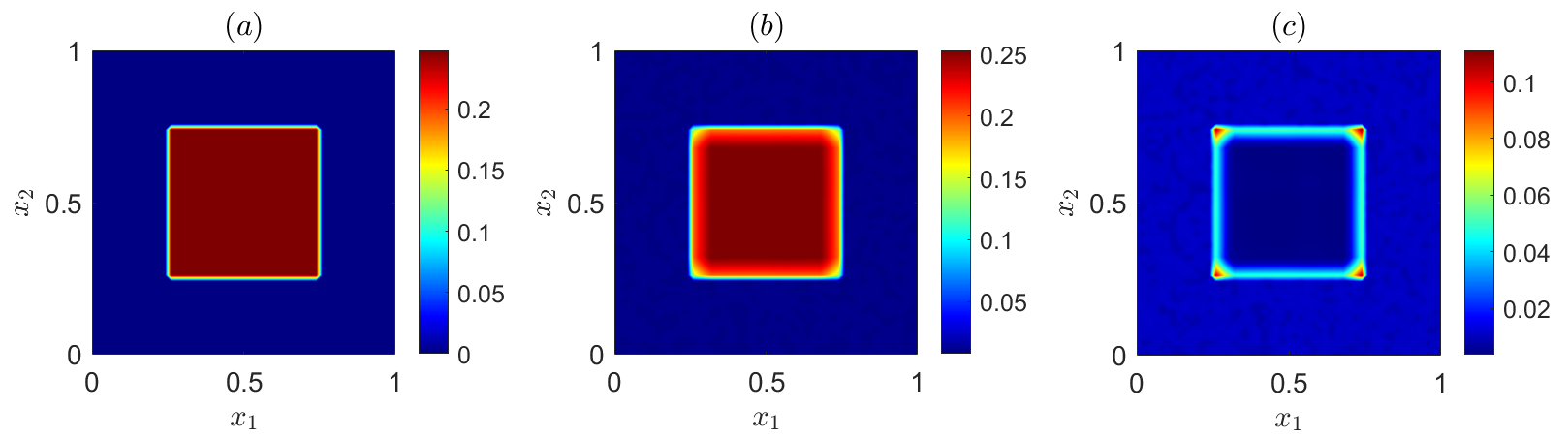}}}
  \centering\caption{\small The computed source functions for Example \ref{examp-2} with $\alpha=0.3$, $\mu(t)=\cos(2\pi t)$, $\omega=[4/50,38/50]^2$ and $\delta_{\mathrm{rel}}=0.5\%$. $(a)$: exact solution $f^*$; $(b)$: regularized solution $f^n$; $(c)$: the absolute error $|f^n-f^*|$.} \label{Fig-2}
\end{figure}

Next we consider two-dimensional problem with region $\overline{\Omega}\times[0,T]=[0,1]^2\times[0,1]$, observation region $\omega=[4/50,38/50]^2$ and relative noise level $\delta_{\mathrm{rel}}=0.5\%$. We fix $\alpha=0.6$ and $\mu(t)=1$,
the exact function we want to recovery is as follows:
\begin{example}\label{examp-2}
$f^*(x_1,x_2)=0.25\chi_{[1/4,3/4]\times[1/4,3/4]}$.
\end{example}
\noindent The space region $\Omega=(0,1)$ is divided into $2\times 40^2$ equal triangles and the time region $[0,1]$ is divided into $50$ equidistant meshes.  For the parameters in Algorithm \ref{alg-acc-pd}, we set $\varsigma=200$, $\upsilon_0=10^{-2}$,
$\gamma =10^{-9}$.
The initial approximation we choose $f^0\equiv 0.1$ and $\rho^0\equiv 0.5$.
Figure \ref{Fig-2} plots the exact solution, regularized solution, and their absolute errors.
The numerical results, as depicted in the figure, exhibit a high level of satisfaction. The recovery solution demonstrates a satisfactory outcome
The precise values concealed within these figures are:
$\delta=6.3437e$--$05$ for the noise level, $n=996$ for the number of iterations, $e_r(f^n,f^*)=0.1187$ for the source error,
and $\mathrm{res}(f^n,g^\delta)=6.9677e$--$05$ for the data error.

\section{Concluding remarks}\label{sec:7}
In this work, we have studied the inverse problem of recovering a space-dependent source term in the TFDEs with the homogeneous Neumann boundary condition.
The approach we adopt in this study differs from most existing literature as we employed an optimal control model with TV regularization, which is beneficial for reconstructing discontinuous or piecewise constant solutions.
By applying the standard Galarkin method with piecewise linear finite element in space and the finite difference scheme in time, a fully discrete scheme for the regularized model was given.
The convergence of the model was analyzed under appropriate assumptions within the framework of fully discrete methods.
Moreover, we reformulated the regularized model as a saddle-point problem and proposed a primal-dual iterative scheme to solve the resulting discretized saddle-point problem.
Finally, the numerical results show that our method is effective.

\bibliographystyle{plain}


\end{document}